\documentclass[12pt]{article}

\usepackage[%
  bookmarks, 
  linkcolor= blue,
  citecolor= blue,
  filecolor= blue,
  urlcolor=  blue,
  colorlinks=true, 
  hyperindex=true,
  pdftitle={Approximation of the Fokker-Planck equation of the stochastic chemostat}, 
  pdfauthor={Fabien Campillo, Marc Joannides and Irène Larramendy-Valverde}
]{hyperref}             

\pdfoutput=1 
\usepackage[utf8x]{inputenc}
\usepackage[english]{babel}
\usepackage{graphicx}           
\usepackage{url}
\usepackage{eurosym}
\usepackage{makeidx}
\usepackage[plain]{algorithm}
\usepackage{algorithmic} 
\usepackage{color} 
\usepackage[margin=3.3cm]{geometry}
\usepackage{setspace}
\usepackage{datetime}
\usepackage{pifont}
\usepackage{tikz}

\usepackage{amsmath,amsfonts,amssymb,latexsym,stmaryrd}

\newcommand{\eqdef}     {\stackrel{{\textrm{\rm\tiny def}}}{=}}

\newtheorem{theorem}      {Theorem}[section]
\newtheorem{theorem*}     {theorem}

\newtheorem{lemma}        [theorem]{Lemma}

\newtheorem{notations}     [theorem]{Notations}
\newcommand{\proof}        {\paragraph{Proof}}

\renewcommand{\Omega}{{\R^2_{+}}}
\newcommand{\mathringOmega}{{\mathring \R^2_{+}}}

\newsavebox{\fmbox}

\newcommand{\N}        {\mathbb N}

\newcommand{\E}        {\mathbb E}

\newcommand{\R}      {\mathbb R}

\renewcommand{\P}      {\mathbb P}

\newcommand{\crochet}[1] {\langle #1 \rangle}

\newcommand{\carre}     {\hfill$\Box$}

\newcommand{\GG}   {{\mathcal G}}

\newcommand{\LL}   {{\mathcal L}}
\newcommand{\NN}   {{\mathcal N}}

\newcommand{\QQ}   {{\mathcal Q}}

\renewcommand{\SS}   {{\mathcal S}}

\newcommand{\rmd}   {{{\textrm{\upshape d}}}}

\renewcommand{\epsilon}{\varepsilon}

\def\dobm{
    \copy1\kern-\wd1\kern0.05ex\copy1\kern-\wd1\kern0.05ex\box1}

\newcommand{\fenumi}  {\textrm{\rm({\textit{i}}\/)}}
\newcommand{\fenumii} {\textrm{\rm({\textit{ii}}\/)}}

\newcommand{\Sin}{s_{\textrm{\tiny in}}}

\newenvironment{psmallmatrix}
{\left(\begin{smallmatrix}}
{\end{smallmatrix}\right)}
\newcommand{\bmax}{b_{\textrm{\tiny\rm max}}}
\newcommand{\smax}{s_{\textrm{\tiny\rm max}}}
\newcommand{\mumax}{\mu_{\textrm{\tiny\rm max}}}
\newcommand{\Sv}{S^{\textrm{\tiny\rm v}}}

\newcommand\blfootnote[1]{%
  \begingroup
  \renewcommand\thefootnote{}\footnote{#1}%
  \addtocounter{footnote}{-1}%
  \endgroup
}

\begin{document}

\title{Approximation of the Fokker-Planck equation \\
of the stochastic chemostat
}

\author{F. Campillo
\thanks{MODEMIC Project-team, INRIA/INRA, UMR MISTEA, Montpellier, France}
\and 
M. Joannides
\thanks{Universit\'e Montpellier 2/I3M, Montpellier, France} 
\footnotemark[1]
\and I. Larramendy-Valverde
\footnotemark[2]
}

\maketitle

\blfootnote{\scriptsize \texttt{Fabien.Campillo@inria.fr}, 
\texttt{marc.joannides@univ-montp2.fr}, 
\texttt{larra@math.univ-montp2.fr}}

\begin{abstract}
We consider a stochastic model of the two-dimensional chemostat as a diffusion process for the concentration of substrate and the concentration of biomass. The model allows for the washout phenomenon: the disappearance of the biomass inside the chemostat. We establish the Fokker-Planck associated with this diffusion process, in particular we describe the boundary conditions that modelize the washout. We propose an adapted finite difference scheme for the approximation of the solution of the Fokker-Planck equation.
\end{abstract}

\section{Introduction}

Many biotechnological processes are modelized with the help of ordinary differential equations (ODE). For example, the dynamic for a single species/single substrate chemostat is classically modelized as \cite{smith1995a}:
\begin{subequations}
\label{eq.x}
\begin{align}
\label{eq.x.s}
  \dot s(t) 
  &= 
  -k\,\mu(s(t))\,b(t)+D\,(\Sin-s(t))\,,
  \\
\label{eq.x.b}
  \dot b(t) 
  &= 
  \{\mu(s(t))-D\}\,b(t)
\end{align}
\end{subequations}
where $b(t)$ and $s(t)$ are the concentrations of biomass and substrate at time $t$
inside the chemostat. The parameters are the dilution rate $D$, the input substrate concentration $\Sin$, and the stoichiometric coefficient $k$.
The specific growth function 
$\mu(s)$ could be of the Monod (non-inhibitory) type:
\begin{align}
\label{eq.monod}
  \mu(s) &= \frac{\mumax\,s}{k_{s}+s} \,,
\end{align}
where $\mumax$ is the maximum growth rate  and $k_{s}$ is the half-saturation; it could also be of the Haldane (inhibitory) type:
\begin{align}
\label{eq.haldane}
  \mu(s) &= \frac{\bar\mu\,s}{k_{s}+s+s^2/\alpha} \,.
\end{align}

\medskip

As pointed out in \cite{campillo2011chemostat}, the system \eqref{eq.x} is simple and applicable to many situations, it can be seen as a limit model of a stochastic birth and death process in high population size asymptotic. Hence \eqref{eq.x} can give account for the mean behavior of the underlying stochastic process but it cannot give account for its the variance. Moreover \eqref{eq.x} fails to propose a realistic representation of the chemostat in small population scenario, that is in cases close to the washout (corresponding to the disappearance of the biomass, i.e. $b(t)=0$).

We present the stochastic model in Section \ref{sec.stochastic.model} and derive the  associated Fokker-Planck equation in Section 
\ref{sec.FP}. A finite difference scheme approximation is detailled in Section  
\ref{sec.DF} and some numerical tests are presented in Section
\ref{sec.numerics}.

\section{The stochastic chemostat model}
\label{sec.stochastic.model}

Consider the stochastic process $X_{t}=(X^1_{t},X^2_{t})=(S_{t},B_{t})$ solution of:
\begin{subequations}
\label{eq.X}
\begin{align}
\label{eq.X.S}
  \rmd S_{t} 
  &= 
  \big\{-k\,\mu(S_{t})\,B_{t}+D\,(\Sin-S_{t})\big\}\,\rmd t 
     + c_{1}\,\sqrt{S_{t}}\,\rmd W^1_{t}\,,
  \\
\label{eq.X.B}
  \rmd B_{t} 
  &= 
  \{\mu(S_{t})-D\}\,B_{t}\,\rmd t + c_{2}\,\sqrt{B_{t}}\,\rmd W^2_{t}\,,
\end{align}
\end{subequations}
where $B_{t}$ and $S_{t}$ are the concentrations of biomass and substrate at time $t$; $W_{t}^1$ and $W_{t}^2$ are independent scalar standard Brownian motions;  $c_{1}>0$ and $c_{2}>0$ are the noise intensities; $W^1_{t}$ and $W^2_{t}$ are independent scalar standard Wiener processes. We suppose that $S_{0}\geq 0$ and $B_{0}\geq 0$ so that $S_{t}\geq 0$ and $B_{t}\geq 0$ for all $t\geq 0$.

\bigskip

The precise analysis of the behavior of the solution of \eqref{eq.X} will be addressed in a forthcoming work \cite{campillo2012analysischemostat}. Still we can describe it simply with some highlights about the classic Cox-Ingersoll-Ross model. Consider the one--dimensional SDE:
\begin{align}
\label{eq.cir}
  \rmd \xi_{t}
  =
  (a+b\,\xi_{t})\,\rmd t
  +
  \sigma\,\sqrt{\xi_{t}}\,\rmd W_{t}
  \,,\quad \xi_{0}=x_{0}\geq 0\,.
\end{align}
with $a\geq 0$, $b\in\R$, $\sigma>0$. According to \cite[Prop. 6.2.4]{lamberton1996a}, for all $x_{0}\geq 0$, $\xi_{t}$ is a continuous process taking values in $\R^{+}$, and let $\tau=\inf\{t\geq 0,\, \xi_{t}=0\}$, then:
\begin{enumerate}
\item
If $a\geq \sigma^2/2$, then $\tau=\infty$ $\P_{x}$--a.s.;
\item
if $0\leq a< \sigma^2/2$ and $b\leq 0$ then $\tau<\infty$ $\P_{x}$--a.s.;
\item
if $0\leq a< \sigma^2/2$ and $b> 0$ then $\P_{x}(\tau<\infty)\in(0,1)$.
\end{enumerate}
In the first case, $\xi_{t}$ never reaches 0. In the second case $\xi_{t}$ a.s. reaches the state 0, in the third case it may reach 0. If $a=0$ then the state 0 is absorbing. 

\medskip

In case of the System \eqref{eq.X},
it is clear that $B=0$ is an absorbing state for \eqref{eq.X.B}, and when $B=0$, \eqref{eq.X.S} reduces to the substrate dynamics conditionally of the washout, namely:
\begin{align}
\label{eq.Sv}
  \rmd \Sv_t
  &=
  D\,(\Sin - \Sv_t)\,\rmd t
  +
  c_{1}\,\sqrt{\Sv_{t}}\,\rmd W^1_t
\end{align}
hence  the solution of this SDE will stay on the half-line
$[0,\infty)$ and:
\begin{enumerate}
\item 
if $D\,\Sin \geq \frac{c_{1}^2}{2}$ then $S_{t}$ never reaches $0$;
\item
if $D\,\Sin < \frac{c_{1}^2}{2}$ then $S_{t}$ reaches $0$ in finite time and is reflected.
\end{enumerate}
Note that, as $c_1$ is ``small'', condition \fenumi\ is more realistic than condition \fenumii: indeed, with a continuous input $\Sin$ , there is no reason for the substrate concentration in the chemostat to vanish.

\bigskip

Simulation schemes for \eqref{eq.X} should respect the previous properties, an adequate choice is:
\begin{subequations}
\label{eq.sim.X}
\begin{align}
\label{eq.sim.X.S}
  S_{t+\delta} 
  &= 
  \big[
  S_{t}
  +
  \big\{-k\,\mu(S_{t})\,B_{t}+D\,(\Sin-S_{t})\big\}\,\delta 
     + c_{1}\,\sqrt{S_{t}}\,\sqrt{\delta}\,w^1_{t}
  \big]_{+}\,,
  \\
\label{eq.sim.X.B}
  B_{t+\delta} 
  &= 
  \big[
  B_{t} 
  +
  \{\mu(S_{t})-D\}\,B_{t}\,\delta + c_{2}\,\sqrt{B_{t}}\,\sqrt{\delta}\,w^2_{t}
  \big]_{+}\,,
\end{align}
\end{subequations}
where $\{w^1_{i\delta}\}_{i\in\N}$ and $\{w^2_{i\delta}\}_{i\in\N}$ are i.i.d. $N(0,1)$ random variables, also independent from $X_{0}$.  Note that $B_{t}=0$ is absorbing for \eqref{eq.sim.X.B}.

\begin{notations}
Let $x=(x_{1},x_{2})=(s,b)\in \R^2_{+}=[0,\infty)^2$ and
\begin{align*}
  f_{1}(x) = f_{1}(s,b) & \eqdef -k\,\mu(s)\,b+D\,(\Sin-s)\,, 
  &
  \sigma_{1}(x)   = \sigma_{1}(s,b)= \sigma_{1}(s) &\eqdef c_{1}\,\sqrt{s}\,, 
\\
  f_{2}(x) = f_{2}(s,b) &\eqdef [\mu(s)-D]\,b\,, 
  &
  \sigma_{2}(x) = \sigma_{2}(s,b)= \sigma_{2}(b) &\eqdef c_{2}\,\sqrt{b}\,, 
\end{align*}
so that  \eqref{eq.X} reads:
\begin{align*}
  \rmd X_{t} &= f(X_{t})\,\rmd t+\sigma(X_{t})\,\rmd W_{t}
\end{align*}
with 
$f(x)=\begin{psmallmatrix} f_{1}(x)  \\  f_{2}(x) \end{psmallmatrix}$,
$\sigma(x)=\begin{psmallmatrix} \sigma_{1}(x) & 0 \\ 0 & \sigma_{2}(x) \end{psmallmatrix}$ and $W_{t}=\begin{psmallmatrix} W^1_{t} \\ W^2_{t} \end{psmallmatrix}$. 

\medskip

Let
$\partial \Omega=\Gamma_{1}\cup\Gamma_{2}$ with 
$\Gamma_{1}=\{(s,b)\in[0,\infty)^2\,;\,b=0\}$
and
$\Gamma_{2}=\{(s,b)\in[0,\infty)^2\,;\,s=0\}$. 


\end{notations}

\section{The Fokker-Planck equation}
\label{sec.FP}

Let $\pi_{t}(\rmd x)=\pi_{t}(\rmd s,\rmd b)$ be the distribution law of of $X_{t}=(S_{t},B_{t})$:
\[
   \pi_{t}(A,B)
   =
   \P(S_{t}\in A\,,\ B_{t}\in B)
\] 
for all Borel sets $A,B$ of $[0,\infty)$.
According to \cite{schuss2010a}, $\pi_{t}(\rmd x)$ of $X_{t}$ can be decomposed as:
\begin{align}
\label{eq.lebesgue}
  \pi_{t}(\rmd x)
  =
  \pi_{t}(\rmd s\times\rmd b)
  =
  \delta_{0}(\rmd b)\,q_{t}(s)\,\rmd s
  +
  p_{t}(s,b)\,\rmd s\,\rmd b
\end{align}
indeed the diffusion process ``lives'' in $\Omega$ but never reaches $\Gamma_{2}$ so the distribution law features only a ``regular'' component $p_{t}(s,b)$ that only charges $\mathringOmega$ and a ``degenerate'' component $q_{t}(s)$ that only charges $\Gamma_{1}$.

As $\pi_{t}$ is a probability distribution we get the normalization property:
\[
   \int_{0}^\infty q_{t}(s)\,\rmd s
   +
   \int_{0}^\infty \int_{0}^\infty p_{t}(s,b)\,\rmd s\,\rmd b
   =
   1\,.
\]
and the washout probability at time $t$ is:
\[
   \P(B_{t}=0)
   =
   \int_{0}^\infty q_{t}(s)\,\rmd s
   =
   1
   -
      \int_{0}^\infty \int_{0}^\infty p_{t}(s,b)\,\rmd s\,\rmd b\,.
\]   

\bigskip

The Fokker-Planck equation in a weak form is:
\begin{align}
\label{eq.FP.weak}
  \frac{\rmd}{\rmd t} 
  \iint_{\R^2_{+}} \pi_{t}(\rmd s,\rmd b)\,\phi(s,b)
  &=
  \iint_{\R^2_{+}} \pi_{t}(\rmd s,\rmd b)\,\LL\phi(s,b)
\end{align}
for all test functions $\phi$, where $\LL$ is the infinitesimal generator defined by:
\begin{align}
\nonumber
  \LL \phi(x)
  &=
  \LL \phi(s,b)
\\
\nonumber
  &\eqdef
  \sum_{i=1}^2 f_i(x)\,\phi'_{x_{i}}(x)
  +
  \frac12\,\sum_{i=1}^2 \sigma^2_{i}(x)\,\phi''_{x_{i}^2}(x)
\\
\label{eq.L}
  &
  =
  \textstyle
  f_{1}(s,b)\,\phi'_{s}(s,b)
  +
  f_{2}(s,b)\,\phi'_{b}(s,b)
  +
  \frac{c_{1}^2}{2}\,s\,\phi''_{s^2}(s,b)
  +
  \frac{c_{2}^2}{2}\,b\,\phi''_{b^2}(s,b)\,.
\end{align}

Using the decomposition \eqref{eq.lebesgue}, the Fokker-Planck equation \eqref{eq.FP.weak} reads:
\begin{align}
\nonumber
 &
 \frac{\rmd}{\rmd t} 
 \Big\{
    \int_{0}^\infty q_{t}(s)\,\phi(s,0)\,\rmd s
    +
    \iint_{\R^2_{+}} p_{t}(s,b)\,\phi(s,b)\,\rmd s\,\rmd b
 \Big\}
 =
 \\
\label{eq.FP.weak2}
 &\qquad\qquad\qquad\qquad
 =
 \int_{0}^\infty q_{t}(s)\,\LL\phi(s,0)\,\rmd s
 +
 \iint_{\R^2_{+}}  p_{t}(s,b)\,\LL\phi(s,b)\,\rmd s\,\rmd b
\end{align}

\begin{lemma}
\label{lemma.LL}
For all functions $\phi\in H^2_{\Gamma_{2}}(\Omega)$ (i.e. $\phi\in H^1(\Omega)$ and $\phi|_{\Gamma_{2}}=0$) and $t\geq 0$
\begin{align*}
 \crochet{p_{t},\LL\phi}
 &
 =
 \int_\Omega \LL^* p_{t}(x)\,\phi(x)\,\rmd x
 +
 \frac{c_{2}^2}{2}\,
 \int_0^\infty p_{t}(s,0)\,\phi(s,0) \,\rmd s
\end{align*}
where $\LL^*$ is the adjoint operator:
\begin{align*} 
\LL^*\psi(x)
 \eqdef
 \textstyle
 - [\psi(x)\,f_{1}(x)]'_{s}
 - [\psi(x)\,f_{2}(x)]'_{b}
 +
 \frac{c_{1}^2}{2}\,[\psi(x)\,s]''_{s^2}
 +
 \frac{c_{2}^2}{2}\,[\psi(x)\,b]''_{b^2}\,.
\end{align*}
\end{lemma}

\proof
By definition of $\LL$:
\begin{align*}
 \crochet{p_{t},\LL\phi}
 =
 \textstyle
 \int_{\Omega} p_{t}(x)\,\LL\phi(x)\,\rmd x
 &=
 \textstyle
 \int_\Omega p_{t}(x)\,f_{1}(x)\,\phi'_{s}(x)\,\rmd x
 +
 \int_\Omega p_{t}(x)\,f_{2}(x)\,\phi'_{b}(x) \,\rmd x
\\
 &\qquad
 \textstyle
 +
 \frac{c_{1}^2}{2}\,
   \int_\Omega p_{t}(x)\,s\,\phi''_{s^2}(x)\,\rmd x
 +
 \frac{c_{2}^2}{2}\,
   \int_\Omega p_{t}(x)\,b\,\phi''_{b^2}(x) \,\rmd x
\end{align*}
we consider separately these four last terms.

From Green’s formula  \cite{brezis2010a}:
 $\int_{\Omega} u'_{x_{i}}\,v\,\rmd x
= -\int_{\Omega}u\,v'_{x_{i}}\,\rmd x + \int_{\partial \Omega} u\,v\,n_{i}\,\rmd\SS_{x}$ where $n_{i}$ is the $i$th component of the outward unit normal $n$,
i.e. $n_{1}(x)=0$ on $\Gamma_{1}$ and $-1$ on $\Gamma_{2}$
and 
$n_{2}(x)=-1$ on $\Gamma_{1}$ and $0$ on $\Gamma_{2}$.
So we get:
\begin{align*} 
 \textstyle
 \int_\Omega p_{t}(x)\,f_{1}(x)\,\phi'_{s}(x)\,\rmd x
 &=
 \textstyle
 -\int_\Omega [p_{t}(x)\,f_{1}(x)]'_{s}\,\phi(x)\,\rmd x
 +\int_{\partial\Omega} p_{t}(x)\,f_{1}(x)\,\phi(x)\,n_{1}(x)\,\rmd \SS_{x}
\\
 &
 \textstyle
 =
 -\int_\Omega [p_{t}(x)\,f_{1}(x)]'_{s}\,\phi(x)\,\rmd x
 -\int_{\Gamma_{2}} p_{t}(x)\,f_{1}(x)\,\phi(x)\,\rmd \SS_{x}
\\
 &
 \textstyle
 =
 -\int_\Omega [p_{t}(x)\,f_{1}(x)]'_{s}\,\phi(x)\,\rmd x\,.
\tag{as $\phi=0$ on $\Gamma_{2}$}
\end{align*}
For the second term:
\begin{align*}
 \textstyle
 \int_\Omega p_{t}(x)\,f_{2}(x)\,\phi'_{b}(x) \,\rmd x
 &=
 \textstyle
 -\int_\Omega [p_{t}(x)\,f_{2}(x)]'_{b}\,\phi(x) \,\rmd x
 +\int_{\partial\Omega} p_{t}(x)\,f_{2}(x)\,\phi(x) \,n_{2}(x)\,\rmd \SS_{x}
\\
 &=
 \textstyle
 -\int_\Omega [p_{t}(x)\,f_{2}(x)]'_{b}\,\phi(x) \,\rmd x
 -\int_{\Gamma_{1}} p_{t}(x)\,f_{2}(x)\,\phi(x)\,\rmd \SS_{x}
\\
 &=
 \textstyle
 -\int_\Omega [p_{t}(x)\,f_{2}(x)]'_{b}\,\phi(x) \,\rmd x\,.
\tag{as $f_{2}=0$ on $\Gamma_{1}$}
\end{align*}
For the third term:
\begin{align*}
\textstyle
   \int_\Omega p_{t}(x)\,s\,\phi''_{s^2}(x)\,\rmd x
 &=
 \textstyle
 -\int_\Omega [p_{t}(x)\,s]'_{s}\,\phi'_{s}(x)\,\rmd x
 +\int_{\partial\Omega} p_{t}(x)\,s\,\phi'_{s}(x)\,n_{1}(x)\,\rmd \SS_{x}
\\
 &=
 \textstyle
 -\int_\Omega [p_{t}(x)\,s]'_{s}\,\phi'_{s}(x)\,\rmd x
 -\int_{\Gamma_{2}} p_{t}(x)\,s\,\phi'_{s}(x)\,\rmd \SS_{x}
\\
 &=
 \textstyle
 - \int_\Omega [p_{t}(x)\,s]'_{s}\,\phi'_{s}(x)\,\rmd x
\tag{as $s=0$ on $\Gamma_{2}$}
\\
 &=
 \textstyle
   \int_\Omega [p_{t}(x)\,s]''_{s^2}\,\phi(x)\,\rmd x
 -
 \int_{\partial\Omega} [p_{t}(x)\,s]'_{s}\,\phi(x)\,n_{1}(x)\,\rmd \SS_{x}
\\
 &=
 \textstyle
   \int_\Omega [p_{t}(x)\,s]''_{s^2}\,\phi(x)\,\rmd x
 +
   \int_{\Gamma_{2}} [p_{t}(x)\,s]'_{s}\,\phi(x)\,\rmd \SS_{x}
\\
 &=
 \textstyle
   \int_\Omega [p_{t}(x)\,s]''_{s^2}\,\phi(x)\,\rmd x\,.
\tag{as $\phi=0$ on $\Gamma_{2}$}
\end{align*}
For the fourth term:
\begin{align*}
\textstyle
   \int_\Omega p_{t}(x)\,b\,\phi''_{b^2}(x) \,\rmd x
 &=
 \textstyle
 - \int_\Omega [p_{t}(x)\,b]'_{b}\,\phi'_{b}(x) \,\rmd x
 + \int_{\partial\Omega} p_{t}(x)\,b\,\phi'_{b}(x)\,n_{2}(x) \,\rmd \SS_{x}
\\
 &=
 \textstyle
 - \int_\Omega [p_{t}(x)\,b]'_{b}\,\phi'_{b}(x) \,\rmd x
 - \int_{\Gamma_{1}} p_{t}(x)\,b\,\phi'_{b}(x) \,\rmd \SS_{x}
\\
 &=
 \textstyle
 - \int_\Omega [p_{t}(x)\,b]'_{b}\,\phi'_{b}(x) \,\rmd x
\tag{as $b=0$ on $\Gamma_{1}$}
\\
 &=
 \textstyle
   \int_\Omega [p_{t}(x)\,b]''_{b^2}\,\phi(x) \,\rmd x
 -
   \int_{\partial \Omega} [p_{t}(x)\,b]'_{b}\,\phi(x)\,n_{2}(x) \,\rmd \SS_{x}
\\
 &=
 \textstyle
   \int_\Omega [p_{t}(x)\,b]''_{b^2}\,\phi(x) \,\rmd x
 +
   \int_{\Gamma_{1}} [p_{t}(x)\,b]'_{b}\,\phi(x)\, \,\rmd \SS_{x}\,.
\end{align*}
Summing up these identities leads to:
\begin{align*}
 \crochet{p_{t},\LL\phi}
 =
 \textstyle
 \crochet{\LL^*\,p_{t},\phi}
 +
 \frac{c_{2}^2}{2}\,
   \int_{\Gamma_{1}} [p_{t}(x)\,b]'_{b}\,\phi(x)\, \,\rmd \SS_{x}
\end{align*}
finally
\begin{align*}
 \textstyle
   \int_{\Gamma_{1}} [p_{t}(x)\,b]'_{b}\,\phi(x)\, \,\rmd \SS_{x}
 &=
 \textstyle
   \int_{\Gamma_{1}} \big\{[p_{t}(x)]'_{b}\,b+p_{t}(x)\big\}\,\phi(x)\, \,\rmd \SS_{x}
\\
 &=
 \textstyle
   \int_{\Gamma_{1}} p_{t}(x)\,\phi(x)\, \,\rmd \SS_{x}
\\
 &=
 \textstyle
   \int_{0}^\infty p_{t}(s,0)\,\phi(s,0)\, \,\rmd s
\end{align*}
proves the lemma.
\carre

\bigskip

According to Lemma \ref{lemma.LL}, \eqref{eq.FP.weak2} becomes:
\begin{align}
\nonumber
 &
 \frac{\rmd}{\rmd t} 
 \Big\{
    \int_{0}^\infty q_{t}(s)\,\phi(s,0)\,\rmd s
    +
    \iint_{\R^2_{+}} p_{t}(s,b)\,\phi(s,b)\,\rmd s\,\rmd b
 \Big\}
 =
 \\
\nonumber
 &
 \qquad\qquad
 =
 \int_{0}^\infty q_{t}(s)\,\LL\phi(s,0)\,\rmd s
 +
 \iint_{\R^2_{+}}  \LL^* p_{t}(s,b)\,\phi(s,b)\,\rmd s\,\rmd b
\\
\label{eq.FP.weak3}
&
 \qquad\qquad\qquad\qquad
 +
 \frac{c_{2}^2}{2}\,
 \int_0^\infty p_{t}(s,0)\,\phi(s,0) \,\rmd s
\end{align}
Let $\phi(s,b)=\varphi(s)\,\psi(b)$ with $\psi(0)=1$, $\psi(b)=0$ for $b>\epsilon$ and $\psi'(0)=\psi''(0)=0$, after letting $\epsilon\to 0$, the previous equation leads to:
\begin{align}
\label{eq.FP.weak3.q}
 &
 \frac{\rmd}{\rmd t} 
    \int_{0}^\infty q_{t}(s)\,\varphi(s)\,\rmd s
 =
 \int_{0}^\infty q_{t}(s)\,\GG\varphi(s)\,\rmd s
 +
 \frac{c_{2}^2}{2}\,
 \int_0^\infty p_{t}(s,0)\,\varphi(s) \,\rmd s
\end{align}
where
\begin{align}
\label{eq.GG}
   \GG\varphi(s)
   =
   \textstyle
   D\,(\Sin-s)\,\varphi'(s)+\frac{c^2_{2}}{2}\,s\,\varphi''(s)
\end{align}
is the infinitesimal generator of the diffusion $S_{t}$ in washout mode, i.e. of the SDE \eqref{eq.Sv}. As 
\eqref{eq.FP.weak3.q} is valid for all test functions $\varphi$, we get the following equation for $q_{t}(s)$:
\begin{subequations}
\label{eq.FP}
\begin{align}
\label{eq.FP.q}
 \frac{\partial}{\partial t} q_{t}(s)
 &=
 \GG^* q_{t}(s)
 +
 \frac{c_{2}^2}{2}\,p_{t}(s,0)
 \,,
 &\forall t\geq 0\,,\ s\in[0,\infty)
\intertext{the equation for $p_{t}(s,v)$ is}
\label{eq.FP.p}
 \frac{\partial}{\partial t} p_{t}(s,v)
 &=
 \LL^* p_{t}(s,v)
 \,,
 &\forall t\geq 0\,,\ (s,v)\in[0,\infty)^2
\end{align}
The initial condition for  \eqref{eq.FP.q} and \eqref{eq.FP.p} are:
\begin{align}
\label{eq.FP.cond.ini}
  q_{t}(s)&=\rho_{\textrm{\tiny v}}(s)\,,
  &
  p_{t}(s,v)&=\rho(s,b)\,.
\end{align}
\end{subequations}
where $\rho_{\textrm{\tiny v}}(s)\,\rmd s\, \delta_{0}(\rmd b)+\rho(s,b)\,\rmd s\,\rmd b$ is the distribution law of $X_{0}=(S_{0},B_{0})$.

The operators are:
\begin{align}
\label{eq.Lv*}
   \GG^*\varphi(s)
   &=
   \textstyle
   -D\,\big[(\Sin-s)\,\varphi(s)\big]'
   +\frac{c^2_{2}}{2}\,\big[s\,\varphi(s)\big]''\,,
\\
\nonumber
  \LL^* \phi(s,v)
  &
  =
  \textstyle
  -\big[ f_{1}(s,b)\,\phi(s,b) \big]'_{s}
  -\big[ f_{2}(s,b)\,\phi(s,b)\big]'_{b}
\\
\label{eq.L*}
  &
  \qquad\qquad\qquad\qquad
  \textstyle
  +
  \frac{c_{1}^2}{2}\,\big[s\,\phi(s,b)]''_{s^2}
  +
  \frac{c_{2}^2}{2}\,\big[b\,\phi(s,b)]''_{b^2}
\end{align}

\medskip

Finally, the Fokker-Planck equation is a system of PDE's: \eqref{eq.FP.p} for $p_{t}(s,v)$ and \eqref{eq.FP.q} for $q_{t}(s)$, the first one is autonomous, and its solution appears as an input for the second PDE.

\section{Approximation}
\label{sec.DF}

Many finite difference schemes and finite element schemes are adapted to space discretization of the system \eqref{eq.FP}. Here we  use the specific finite difference scheme proposed in  \cite{kushner1977a}. This classical scheme presents nice numerical  properties and it also can be interpreted as an approximation of the solution of 
\eqref{eq.X} by a pure jump Markov process on a finite discretization grid, the resulting system in discrete-space and continuous-time is the exact Fokker-Planck equation (forward Kolmogorov equation) associated with this pure jump process. The infinitesimal generator $\LL$ of the SDE \eqref{eq.X} is given by \eqref{eq.L}, this operator fully characterizes the distribution law of the process $X_{t}=(S_{t},B_{t})$, indeed the set of equations \eqref{eq.FP} is totally determined by the operator $\LL$ as $\GG$ is only the restriction of $\LL$ to $\Gamma_{2}$.

\medskip

\begin{figure}
\begin{center}
\includegraphics[width=8.5cm,keepaspectratio]{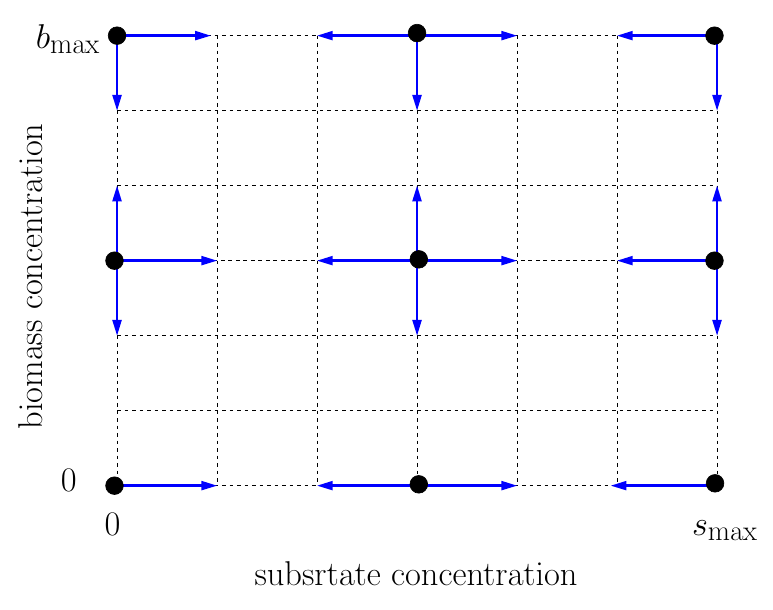}
\end{center}
\caption{Discretized domain $G_{h}$.}
\label{fig.grid}
\end{figure}

The finite difference scheme is detailed in  \ref{sec.DF.appendix}, it leads to the following approximation of the infinitesimal generator:
\[
  \LL \phi(x)
  \simeq
  \LL_{h} \phi(x)
  =
  \sum_{y\in G_{h}}\LL_h(x,y)\,\phi(y)
\]
for $x\in G_{h}$ where:
\begin{align*}
 G_{h}
 &\eqdef 
 \{x=(k_{1}\,h_{1},k_{2}\,h_{2})\,;\,k_{i}=0,\dots,N_{i},\, i=1,2\}
 \,,
\\
 \mathring G_{h}
 &\eqdef 
 \{x=(k_{1}\,h_{1},k_{2}\,h_{2})\,;\,k_{i}=1,\dots,N_{i}-1,\, i=1,2\}
 \,,
\\
 G^1_{h}
 &\eqdef 
 \{x=(k_{1}\,h_{1},0)\,;\,k_{1}=0,\dots,N_{1}\}
 \,,
\end{align*}
are the grid version of $\Omega$, $\mathring \Omega$ and $\Gamma_{1}$ respectively, see Figure \ref{fig.grid}.

\medskip

For the interior points $x\in \mathring G_{h}$ the finite difference scheme is:
\begin{align*}
\left\{
\begin{array}{rll}
  \LL_{h} (x,x)
  &=
  - \frac{|f_1(x)|}{h_1} - \frac{|f_2(x)|}{h_2} - \frac{\sigma^2_1(x)}{h_1^2}
  - \frac{\sigma^2_2(x)}{h_2^2}
  \,,
\\
  \LL_{h} (x,x\pm h_i\,e_i)
  &=
  \frac{f_{i}^\pm(x)}{h_{i}}
  +
  \frac{\sigma^2_{i}(x)}{2\,h_{i}^2} 
     \,,\quad i=1,2\,,
\\[0.7em]
  \LL_{h} (x,y)
  &=
  0 \qquad\textrm{otherwise.}
\end{array}
\right.
\end{align*}
For the boundary points  $x\in G_{h}\setminus \mathring G_{h}$ the finite difference schemes are detailed in \ref{sec.bc}. They correspond to the Figure \ref{fig.grid}: for $s=\smax$ or $b=\bmax$, we must impose reflecting conditions, for for $s=0$ or $b=0$, the boundary conditions are natural, they derive from the value of the coefficients. Indeed, when $b=0$, then $f_{2}=\sigma_{2}=0$ and the jump process stays on the boundary ``$b=0$'' (it cannot jump to $b=h_{2}$ or to $b=-h_{2}$). When $s=0$, then $f_{1}=D\,\Sin$ and $\sigma_{1}=0$, so the jump process can only jump to $s=h_{1}$.
 
\medskip

We obtain a matrix $\LL_{h}=[\LL_{h}(x,y)]_{x,y\in G_{h}}$ which is the infinitesimal generator of a pure jump Markov process $(X^h_{t})_{t\geq 0}$ in continuous time and discrete state space $G_{h}$. Starting from a point $x$ of the grid, the  process $X^h_{t}$ stays there during a time exponentially distributed with parameter $-\LL_{h}(x,x)$ then it jumps to a point $y$ with probability $\LL_{h}(x,y)/(-\LL_{h}(x,x))$ for all $y\in\GG_{h}$, and $\LL_{h}$ is a $\QQ$-matrix as $\sum_{y\in\GG_{h}}\LL_{h}(x,y)=0$. Then the following  Kolmogorov forward equation:
\begin{align}
\label{eq.FP.h}
  \frac{\partial}{\partial t}p^h_{t}(x) 
  &= 
  \LL_{h}^* \, p^h_{t}(x)
\end{align}
 gives the evolution of the distribution law $p^h_{t}$ of $X^h_{t}$, $p^h_{t}(x)=\P(X^h_{t}=x)$, $x\in G_{h}$.
 
\medskip
 
It is important to note that this approach gives an approximation of the coupled system of PDEs \eqref{eq.FP}: 
$(p^h_{t}(x))_{x\in G_{h}\setminus G^1_{h}}$ is an approximation of $(p_{t}(s,b))_{(s,b)\in(0,\infty)^2}$
and
$(p^h_{t}(x))_{x\in G^1_{h}}$ is an approximation of $(q_{t}(s))_{s\in(0,\infty)}$.

\medskip

For the time-discretization we use the
implicit Euler scheme approximation:
\begin{align*}
   \frac{p^h_{t+\delta}(x)-p^h_{t}(x)}{\delta}
   =
   \LL_{h}^*\,p^h_{t+\delta}(x) 
\end{align*}
that is:
\begin{align*}
   (I-\delta\,\LL^*_{h})\,p_{t+\delta}(x)
   =
   p_{t}(x) \,.
\end{align*}

\section{Numerical results}
\label{sec.numerics}

\subsection{Comparison}

Many works \cite{imhof2005a} propose the following structure for the diffusion coefficients:
\begin{subequations}
\label{eq.X.comp}
\begin{align}
\label{eq.X.S.comp}
  \rmd S_{t} 
  &= 
  \big\{-k\,\mu(S_{t})\,B_{t}+D\,(\Sin-S_{t})\big\}\,\rmd t 
     + c_{1}\,S_{t}\,\rmd W^1_{t}\,,
  \\
\label{eq.X.B.comp}
  \rmd B_{t} 
  &= 
  \big\{\mu(S_{t})\,B_{t}-D\,B_{t}\big\}\,\rmd t + c_{2}\,B_{t}\,\rmd W^2_{t}\,.
\end{align}
\end{subequations}
It is slightly different from \eqref{eq.X}. In large population size, these two models are rather equivalent; they differ drastically in the washout regime.

\begin{figure}
\begin{center}
\begin{tabular}{rcc}
\rotatebox{90}{\hskip 3em \scriptsize $t=1$}
&
\includegraphics[trim = 8mm 26mm 12mm 29mm,width=6.5cm,keepaspectratio]
       {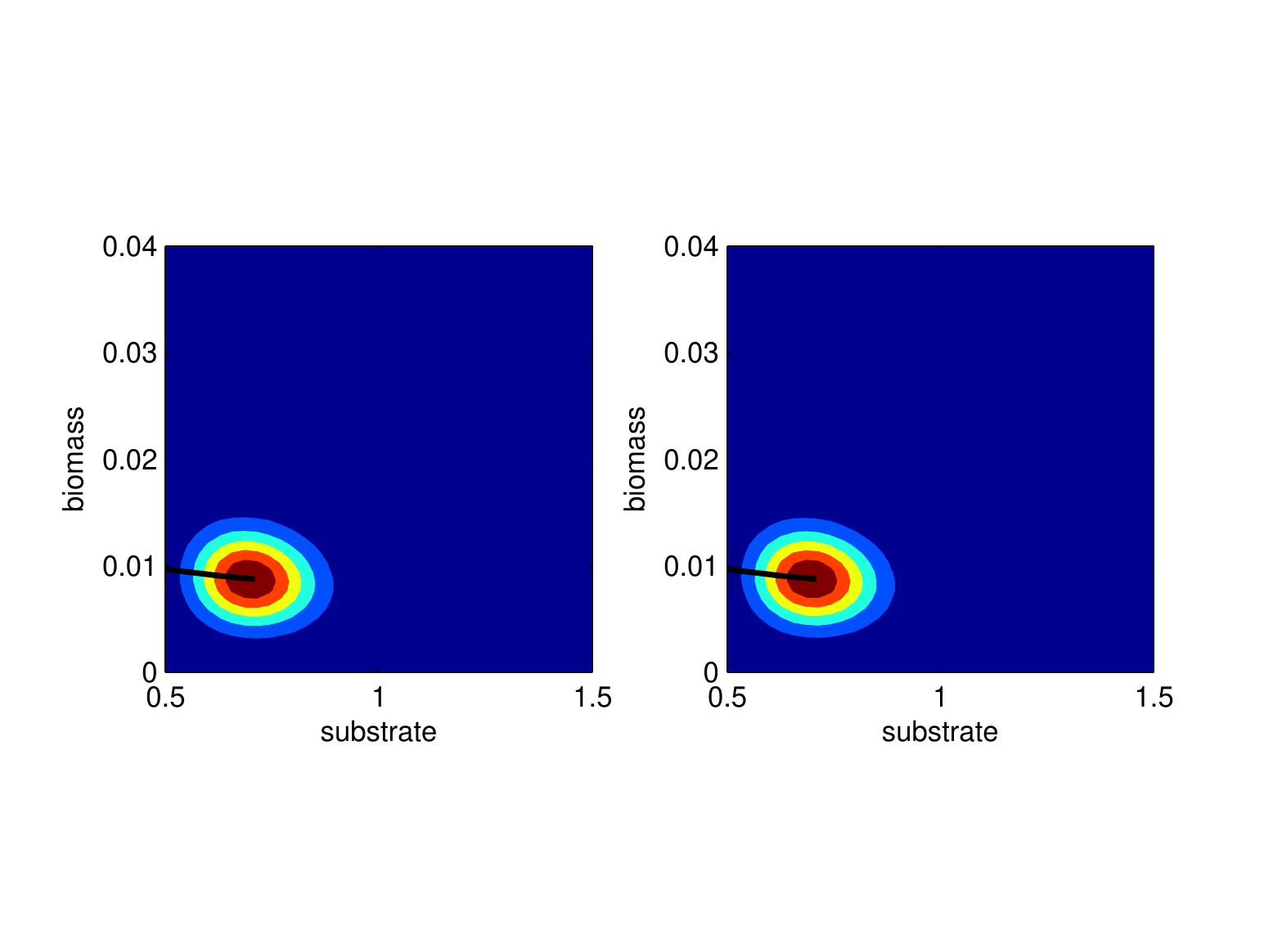}
&
\includegraphics[trim = 8mm 26mm 12mm 29mm,width=6.5cm,keepaspectratio]
       {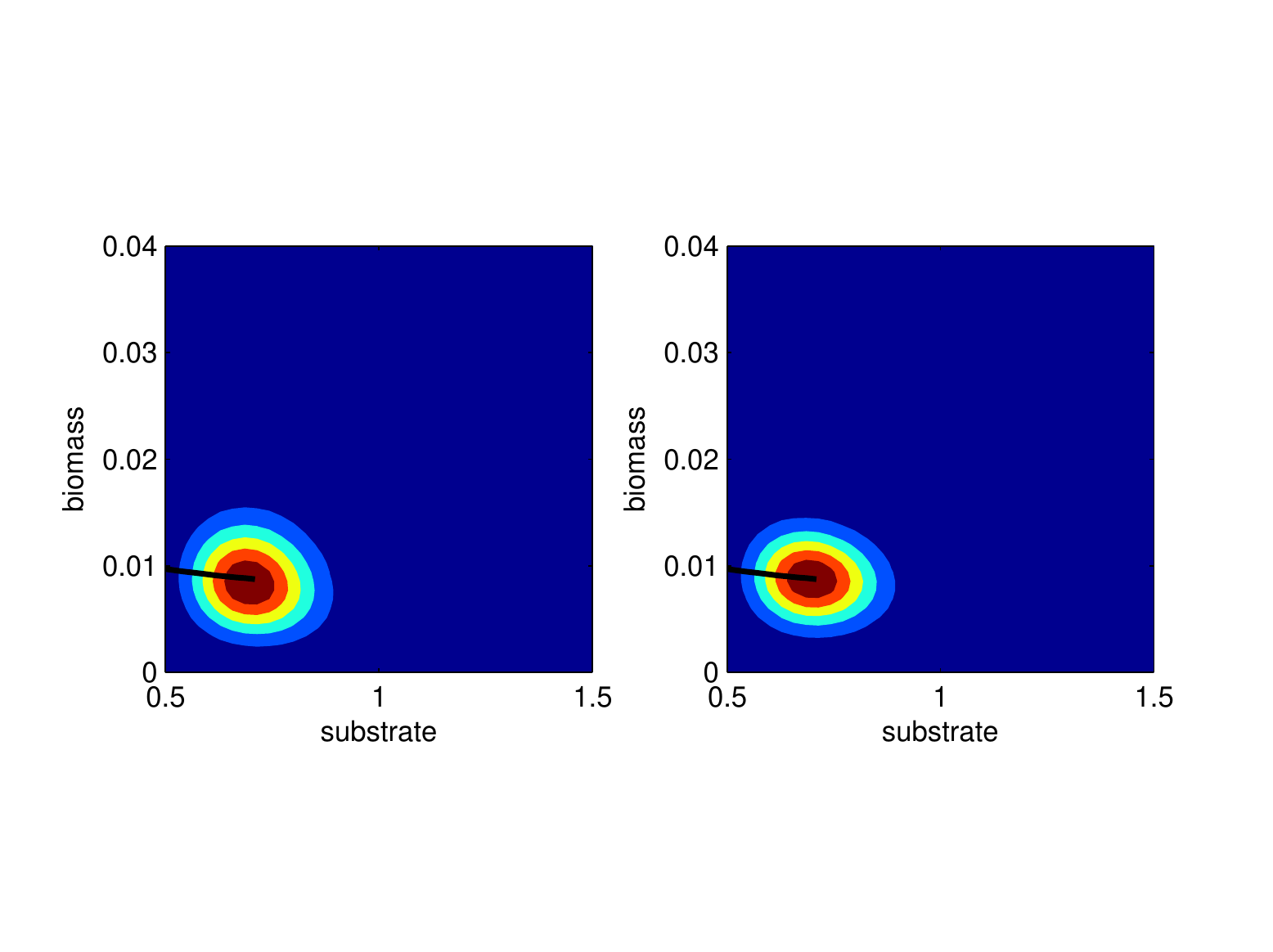}
\\
\rotatebox{90}{\hskip 3em \scriptsize $t=5$}
&
\includegraphics[trim = 8mm 26mm 12mm 29mm,width=6.5cm,keepaspectratio]
       {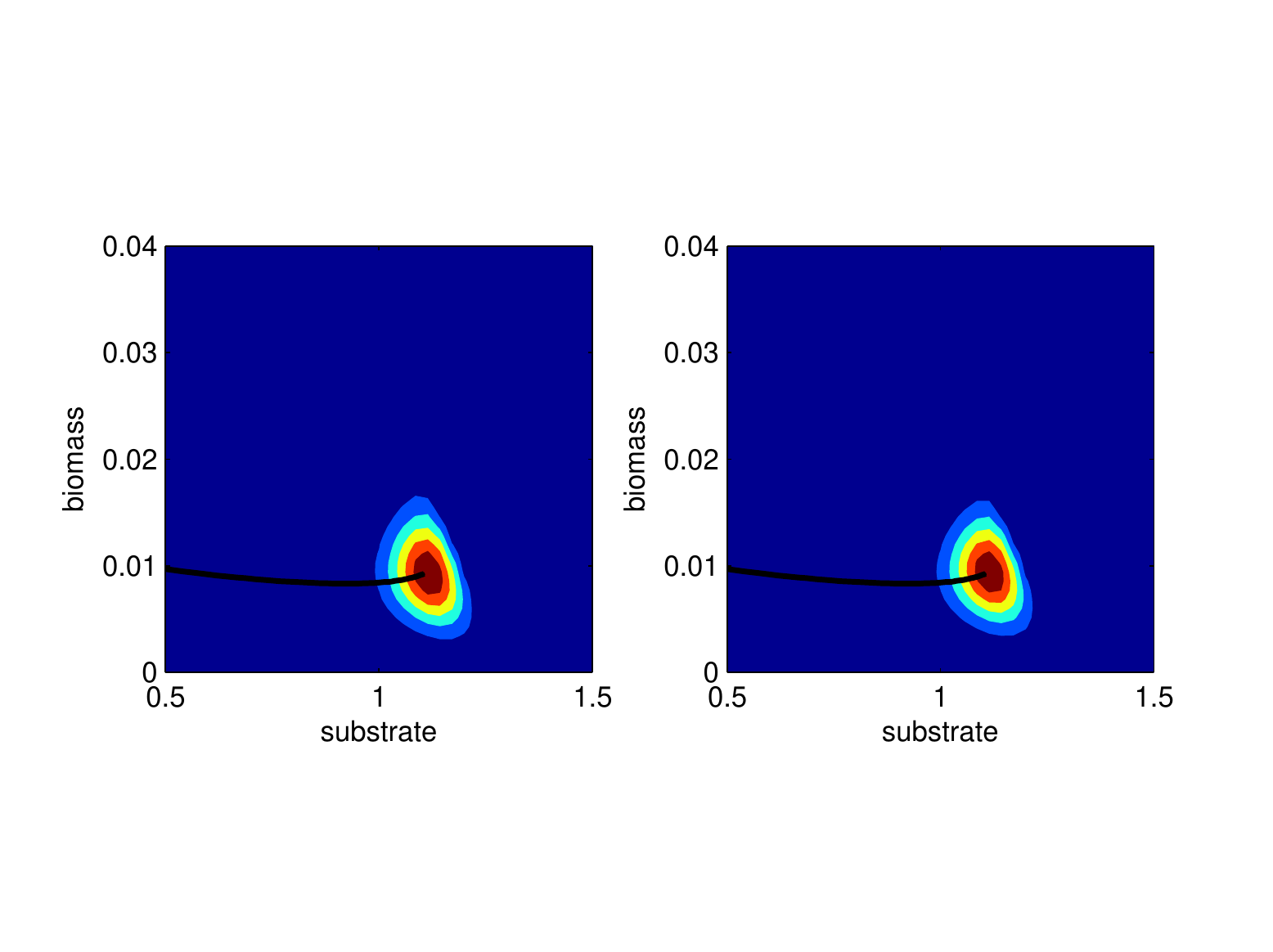}
&
\includegraphics[trim = 8mm 26mm 12mm 29mm,width=6.5cm,keepaspectratio]
       {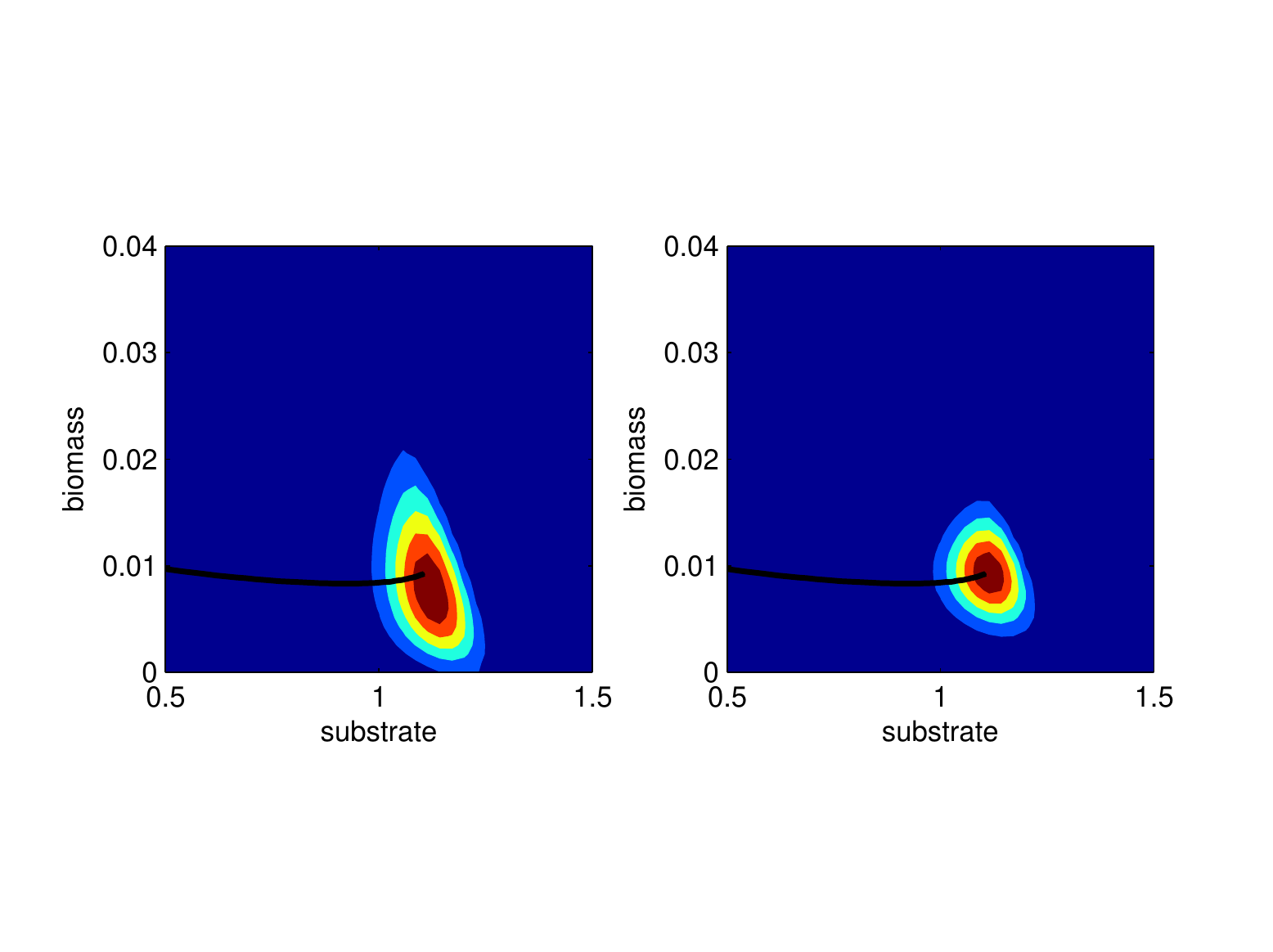}
\\
\rotatebox{90}{\hskip 3em \scriptsize $t=10$}
&
\includegraphics[trim = 8mm 26mm 12mm 29mm,width=6.5cm,keepaspectratio]
       {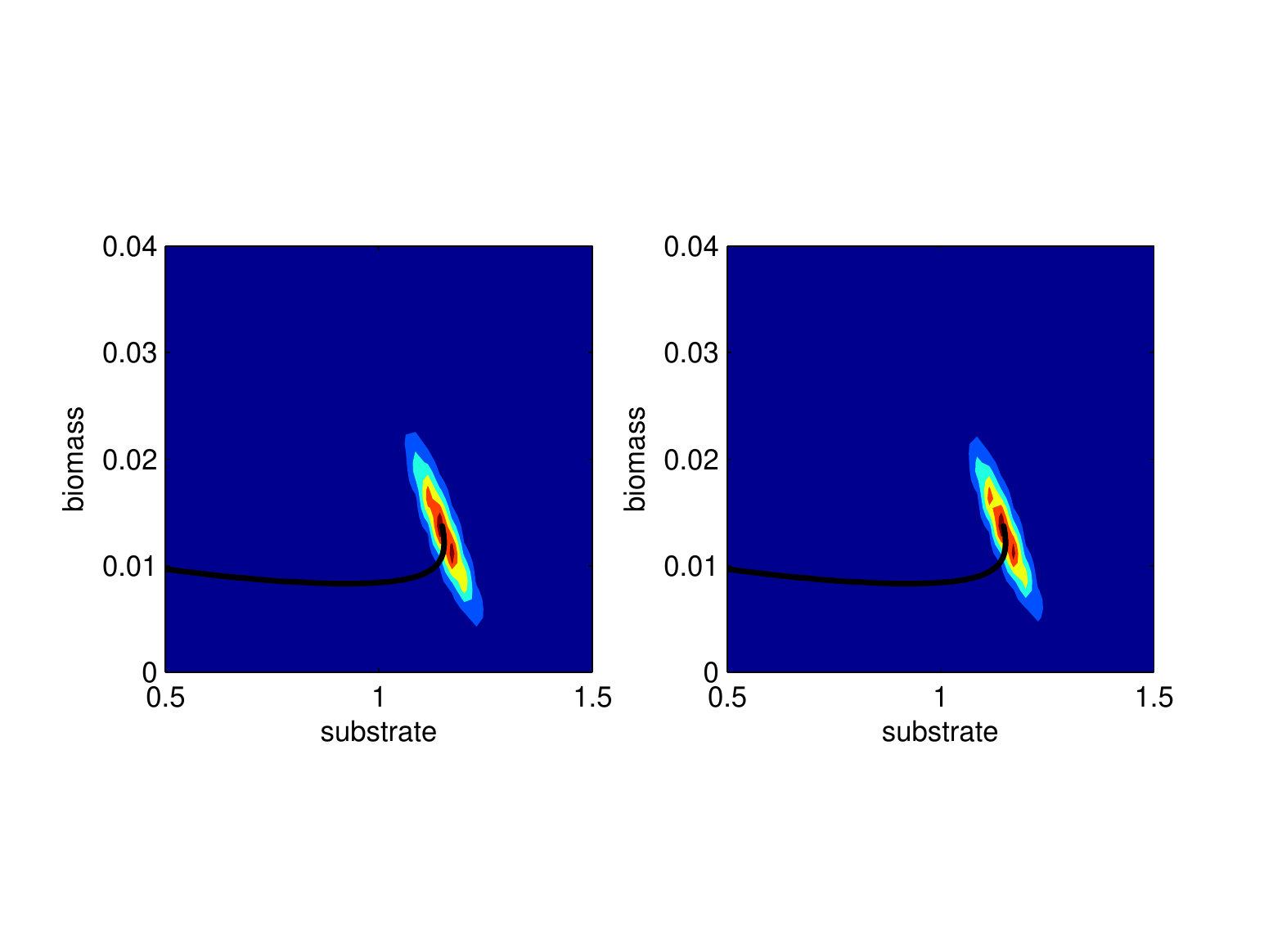}
&
\includegraphics[trim = 8mm 26mm 12mm 29mm,width=6.5cm,keepaspectratio]
       {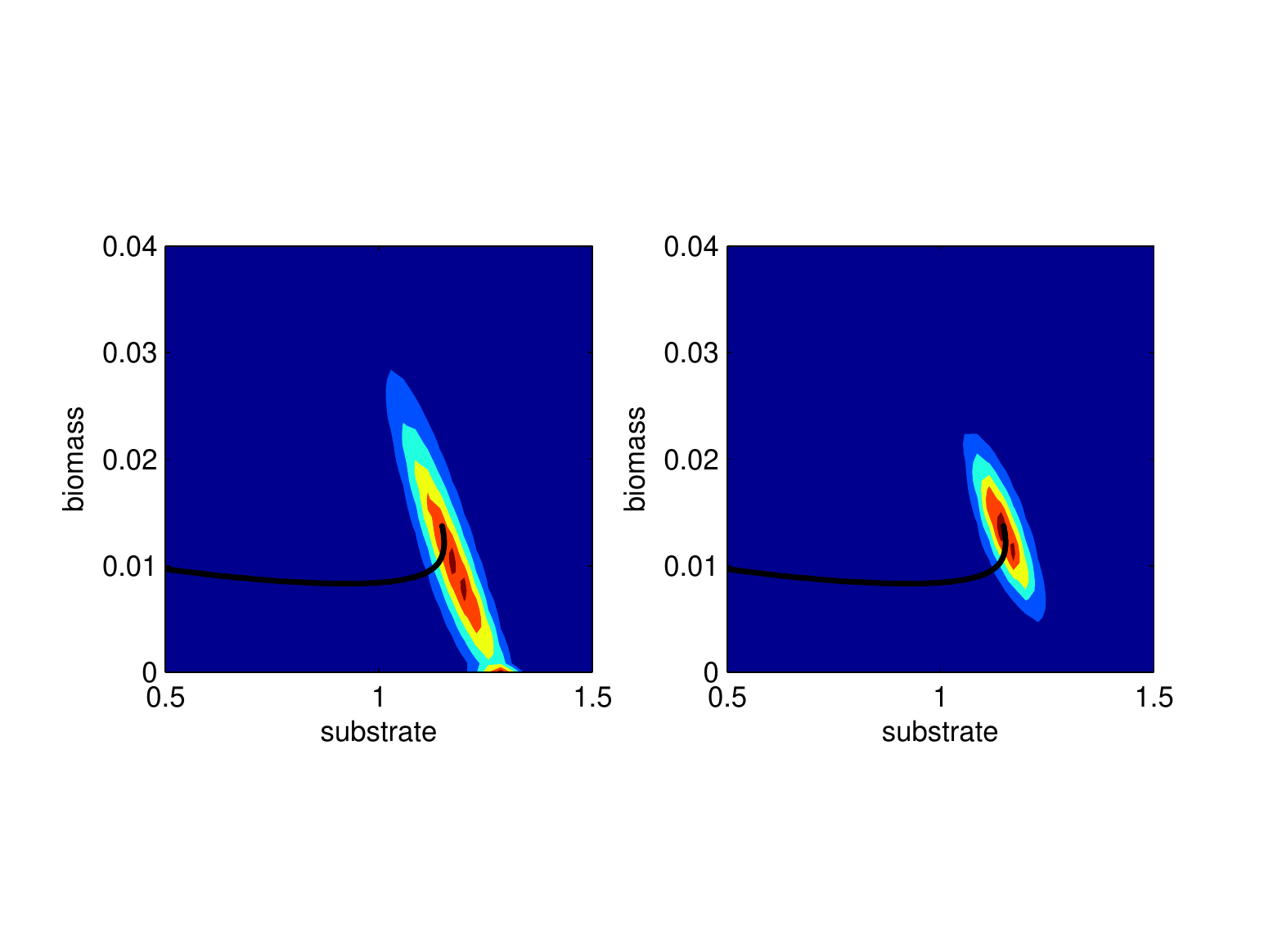}
\\
\rotatebox{90}{\hskip 3em \scriptsize $t=15$}
&
\includegraphics[trim = 8mm 26mm 12mm 29mm,width=6.5cm,keepaspectratio]
       {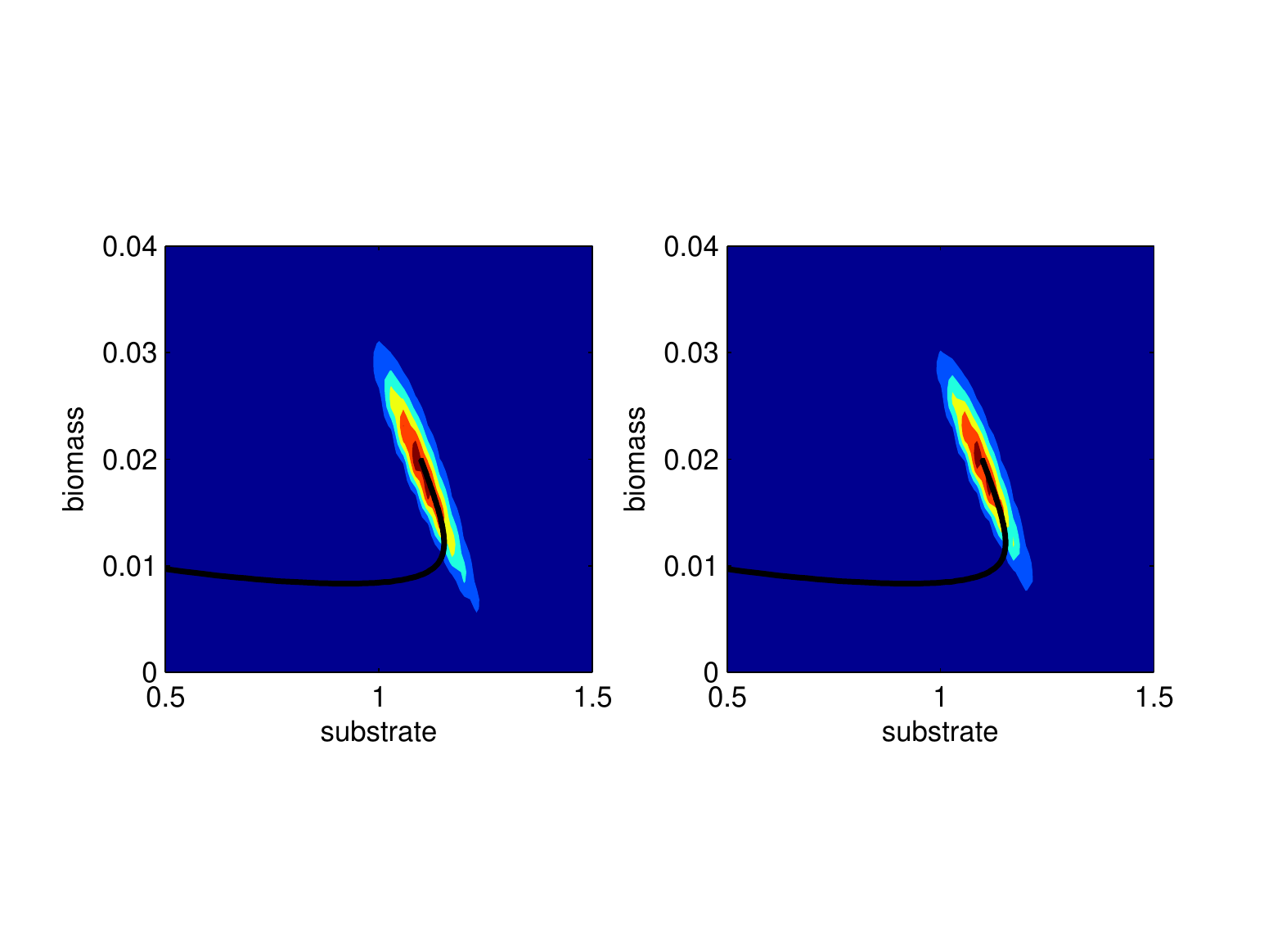}
&
\includegraphics[trim = 8mm 26mm 12mm 29mm,width=6.5cm,keepaspectratio]
       {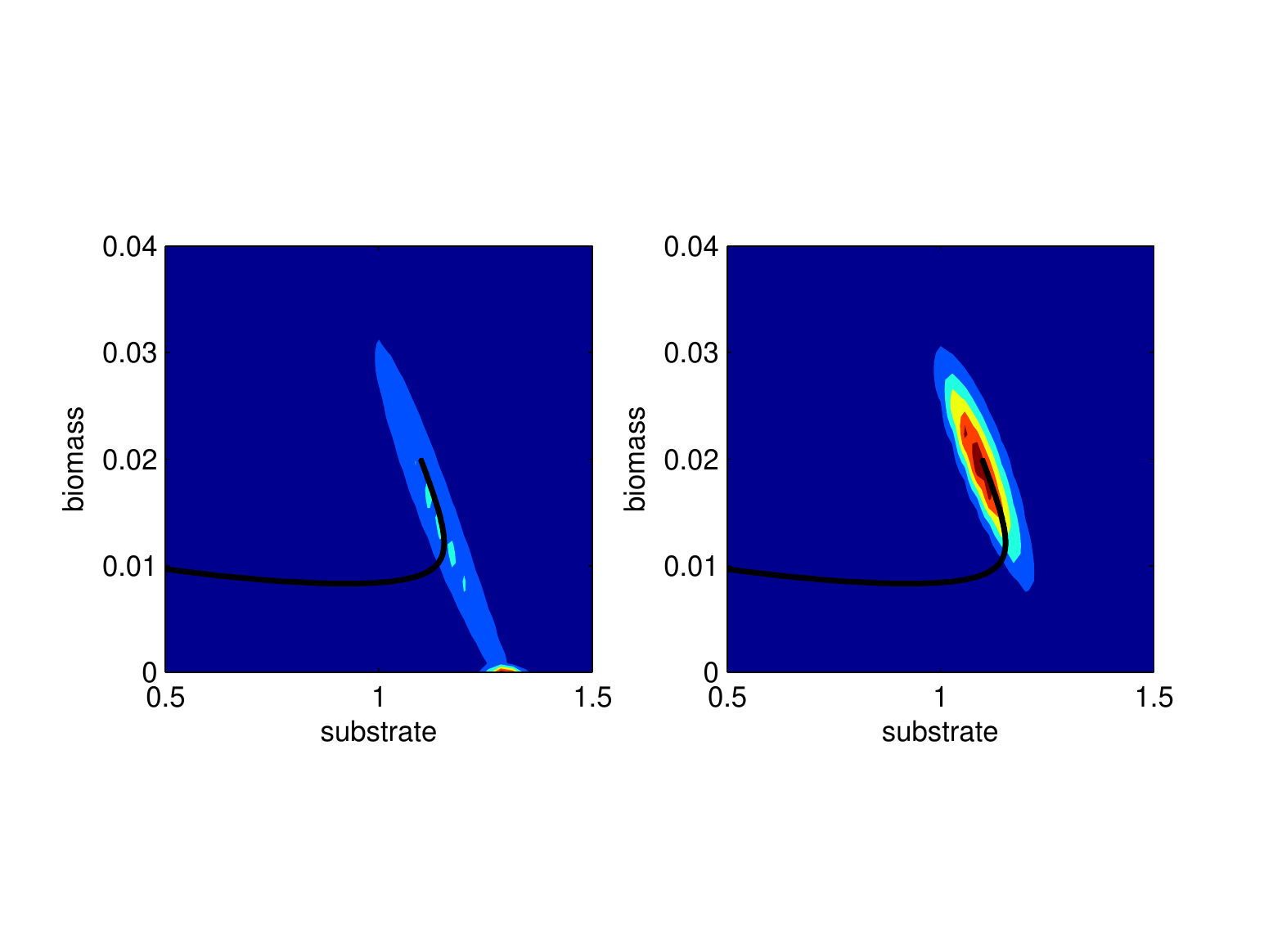}
\\
\rotatebox{90}{\hskip 3em \scriptsize $t=20$}
&
\includegraphics[trim = 8mm 26mm 12mm 29mm,width=6.5cm,keepaspectratio]
       {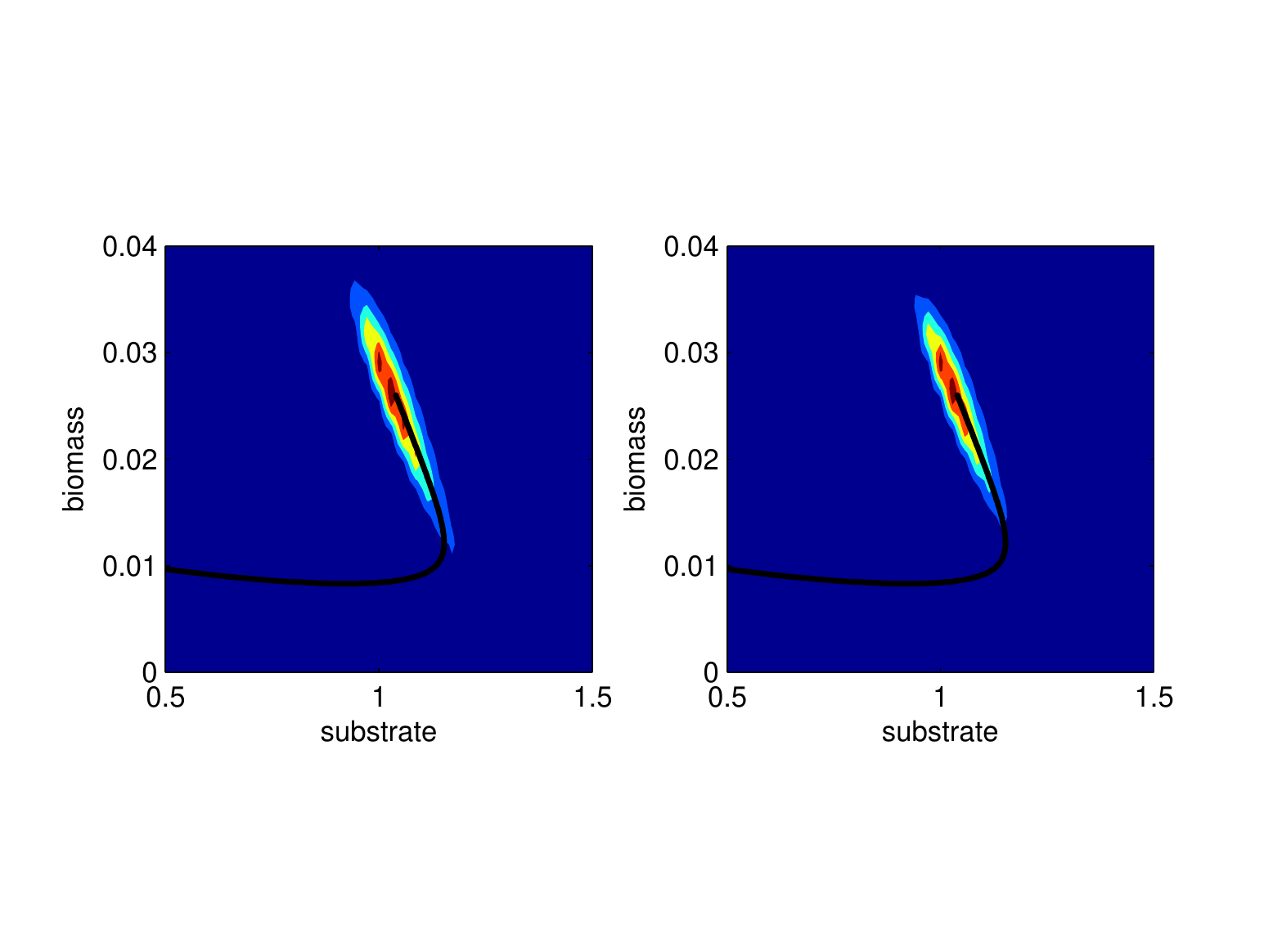}
&
\includegraphics[trim = 8mm 26mm 12mm 29mm,width=6.5cm,keepaspectratio]
       {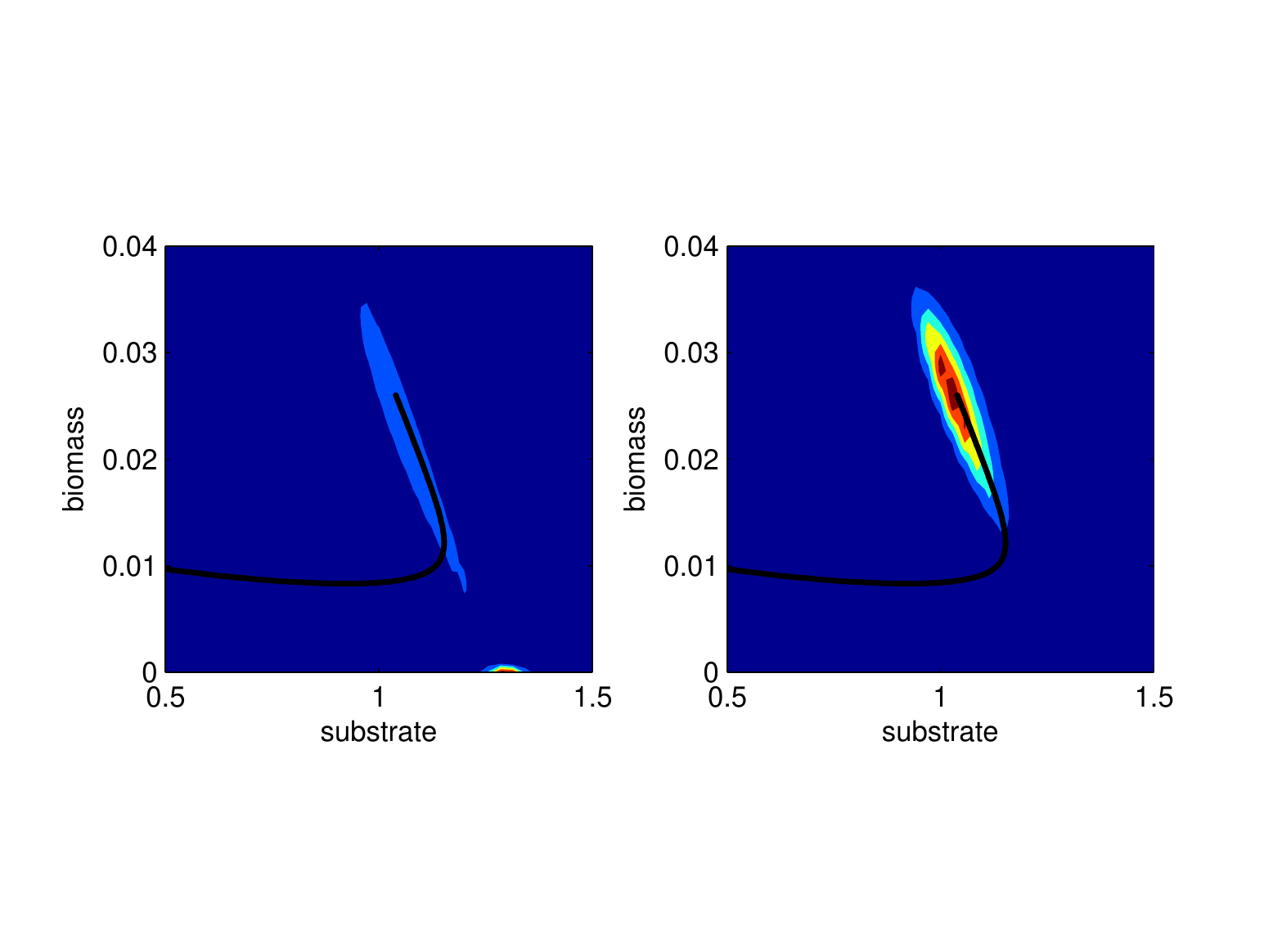}
       
\\
&
{\it\hspace{1em}case 1.a\hspace{6em}case 2.a}
&
{\it\hspace{1em}case 1.b\hspace{6em}case 2.b}
\end{tabular}
\end{center}
\caption{In cases ``1'' the diffusion coefficients are $\sigma_{1}(s)=c_{1}\,\sqrt{s}$ and $\sigma_{2}(b)=c_{2}\,\sqrt{b}$;
in cases ``2'' the diffusion coefficients are $\sigma_{1}(s)=c_{1}\,s$ and $\sigma_{2}(b)=c_{2}\,b$. In cases ``a'' $c_{1}=c_{2}=0.005$; in cases ``b'' $c_{1}=c_{2}=0.02$. For small noise intensities (cases ``a''), cases ``1'' and ``2'' behave rather similarly. For higher noise intensities (cases ``b''), as the law $\pi_{t}$ of $(S_{t},B_{t})$ is closer to the absorbing ``washout'' boundary $\{(s,b)\in\R^2_{+};b=0\}$, cases ``1'' and ``2'' behave rather similarly. See Figure \ref{fig.compare.proba} for the evaluation of the washout probability.}
\label{fig.compare}
\end{figure}

\begin{figure}
\begin{center}
\includegraphics[width=6.5cm,keepaspectratio]
       {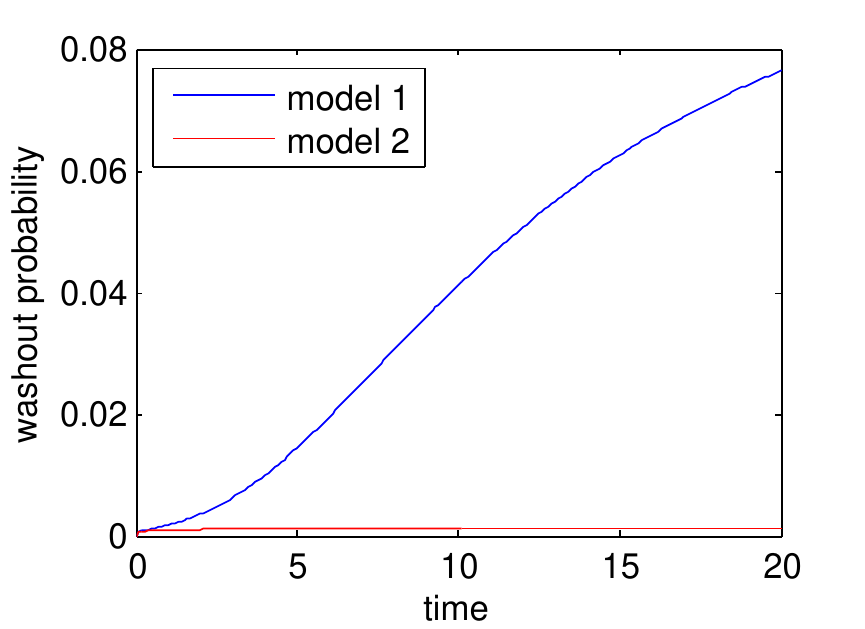}
\end{center}
\caption{Washout probability --- Following Figure \ref{fig.compare}: we compute $t\to \P(B_{t}=0)$ for 
the case ``1'' (model 1: $\sigma_{1}(s)=c_{1}\,\sqrt{s}$, $\sigma_{2}(b)=c_{2}\,\sqrt{b}$) and for case  ``2'' (model 2: $\sigma_{1}(s)=c_{1}\,s$, $\sigma_{2}(b)=c_{2}\,b$).}
\label{fig.compare.proba}
\end{figure}

In this test we use the Monod growth rate function \eqref{eq.monod} and the parameters:
$k=10$, $\Sin=1.3$ (mg/l), $D=0.4$ (1/h), $\mumax=3$ (1/h) , $k_{s}=6$ (mg/l). The initial law is $(S_{0},B_{0}) \sim \NN(0.45,10^{-5})\otimes \NN(0.01,10^{-5})$. The discretization parameters are $\smax=2$, $\bmax=0.06$, $\delta=0.1$, $N_{1}=N_{2}=70$. In Figure \ref{fig.compare}, we see that with small noise intensities the simulation of the two models are very similar; with higher small noise intensities, the simulations are very different. This is due to the fact that the behavior of the two diffusion processes near the boundary ``$b=0$'' are different: with the model \eqref{eq.X} the washout regime is attainable which is not the case with the model 
\eqref{eq.X.comp}. In Figure \ref{fig.compare.proba} we compare the evolution of the washout probability $t\to\P(B_{t}=0)$ for both models, we clearly see that the model 
\eqref{eq.X.comp} does not give account for this probability.

\subsection{Simulation with the Haldane growth rate function}

In this test we use the Haldane growth rate function \eqref{eq.haldane} and the parameters:
$k=2$, $\Sin=2.4$ (mg/l), $D=0.1$ (1/h), $\bar\mu=5$ (1/h) , $k_{s}=10$ (mg/l), $\alpha=0.03$: $c_{1}=c_{2}=0.01$. The initial law is $(S_{0},B_{0}) \sim \NN(1.5,10^{-5})\otimes \NN(0.68,10^{-5})$. The discretization parameters are $\smax=3$, $\bmax=2.5$, $\delta=0.25$, $N_{1}=N_{2}=300$.

\begin{figure}
\begin{center}
\includegraphics[width=6.5cm,keepaspectratio]
       {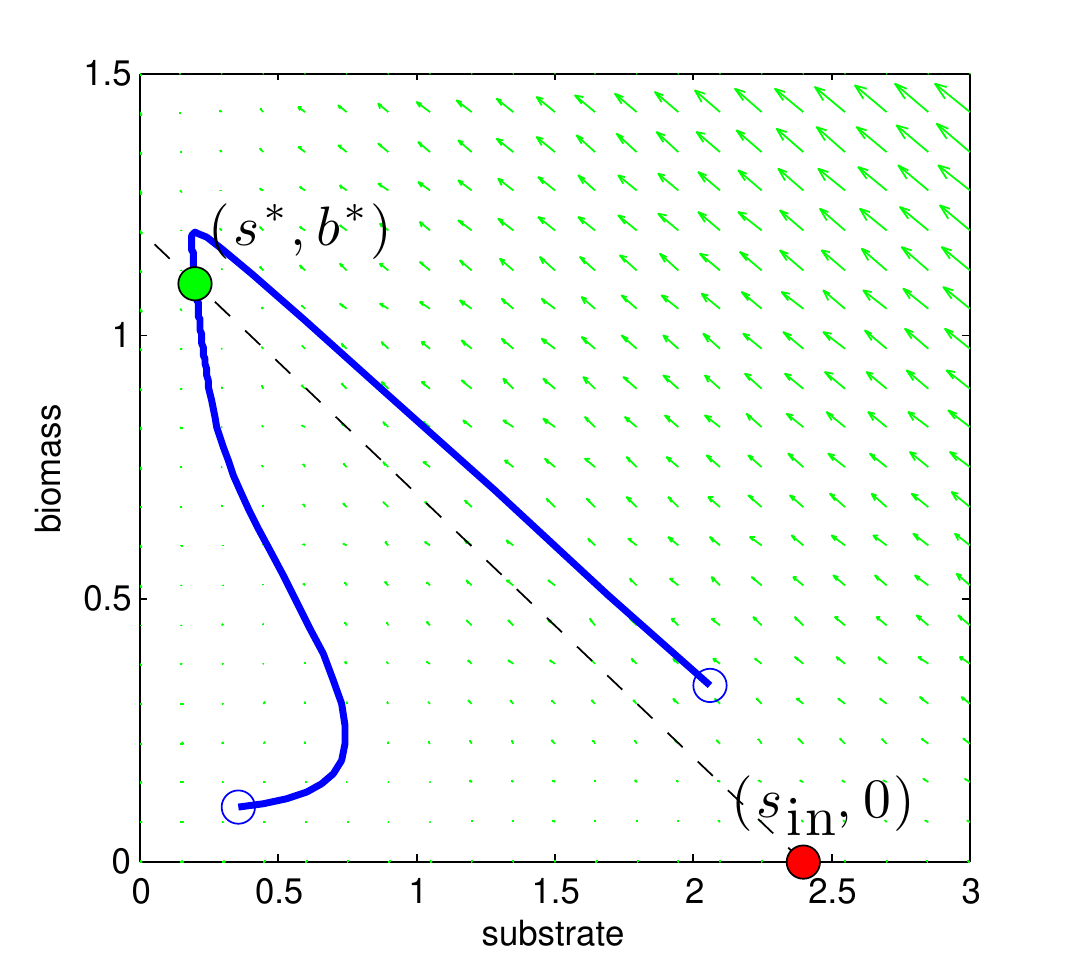}
\includegraphics[width=6.5cm,keepaspectratio]
       {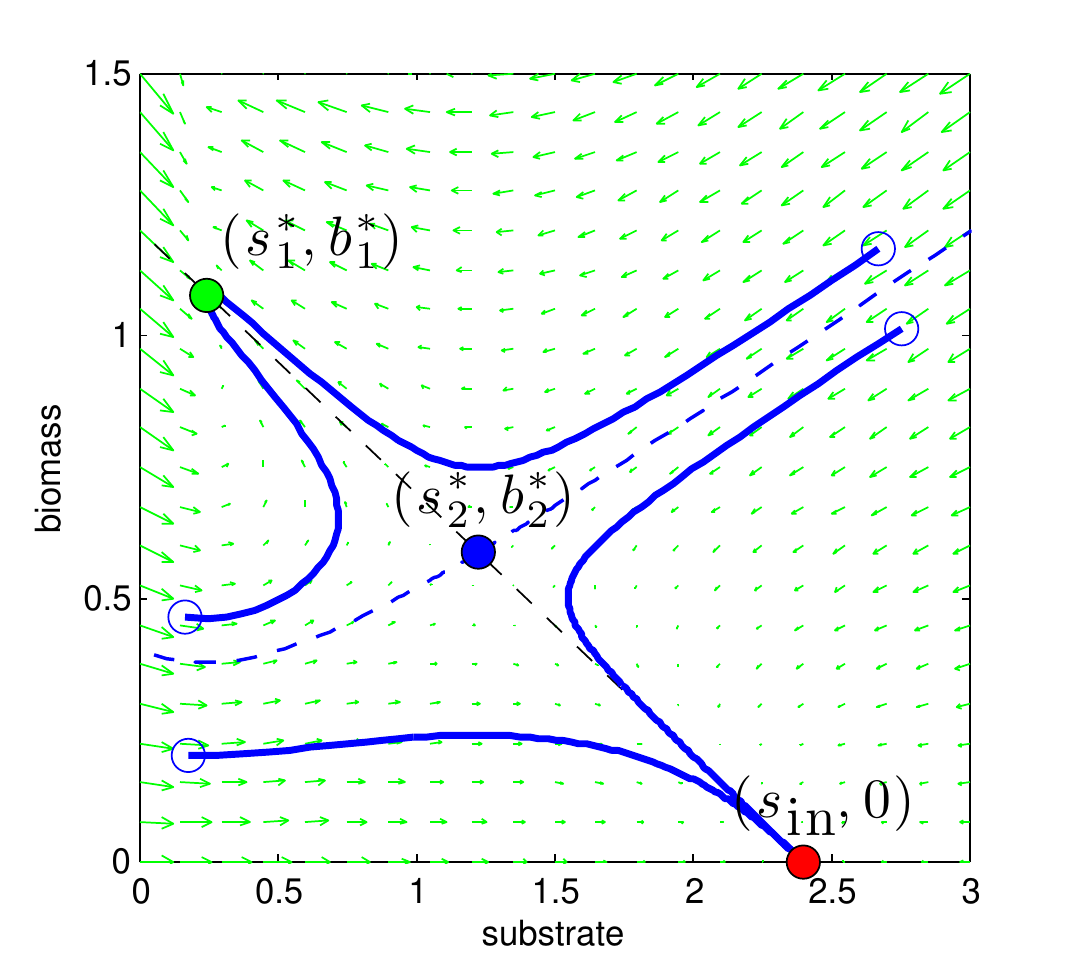}
\end{center}
\caption{Phase portraits for the system \eqref{eq.x} for the Monod growth function (left) and the Haldane growth function. {\bf Left (Monod case)}:  there are two equilibrium states: the washout equilibrium (red dot) is unattractive, the equilibrium point $(s^*,b^*)$ with 
$s^*=k_s\,D/(\mumax-D)$ (solution of $\mu(s)=D$) and 
$b^*=(\Sin-s^*)/k$ is attractive. We suppose that $\mumax>D$. The dashed line is $ b=(\Sin-s)/k$, in blue two trajectories (blue circles: initial positions).
{\bf  Right (Haldane case)}: the washout is still an equilibrium point but now it is attractive, there are two other equilibrium points given as solutions of $\mu(s)=D$ (we suppose that it admits two separate solutions), $(s^*_{1},b^*_{1})$ is attractive (corresponding to the smallest value of $s$), $(s^*_{2},b^*_{2})$ is unattractive. The black dashed curve separates the two basins of attraction; in blue four trajectories (blue circles: initial positions).}
\label{fig.monod.haldane}
\end{figure}

\begin{figure}
{\small
\begin{center}
\begin{tabular}{ccc}
\includegraphics[width=4.7cm,keepaspectratio]{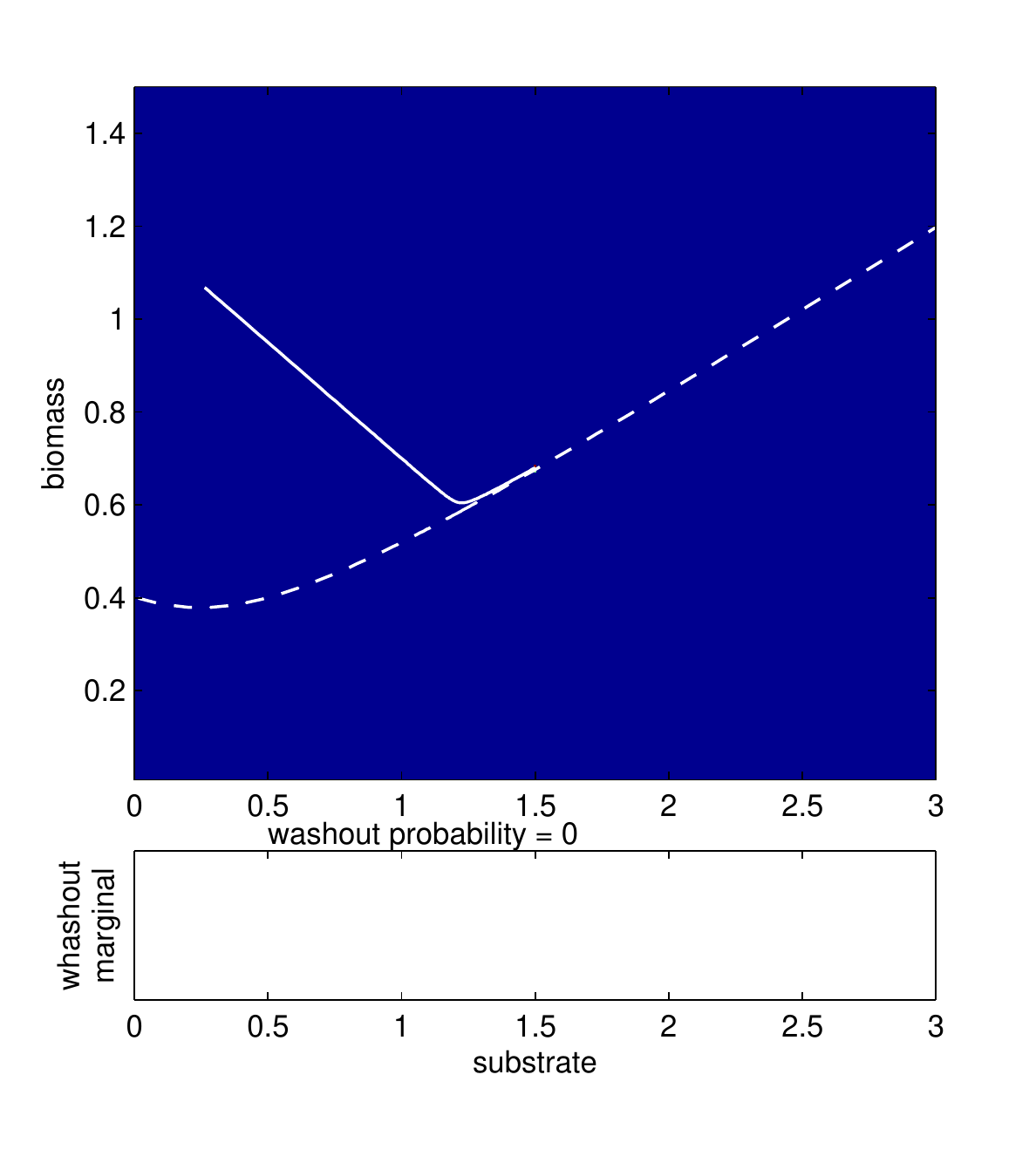}
&
\includegraphics[width=4.7cm,keepaspectratio]{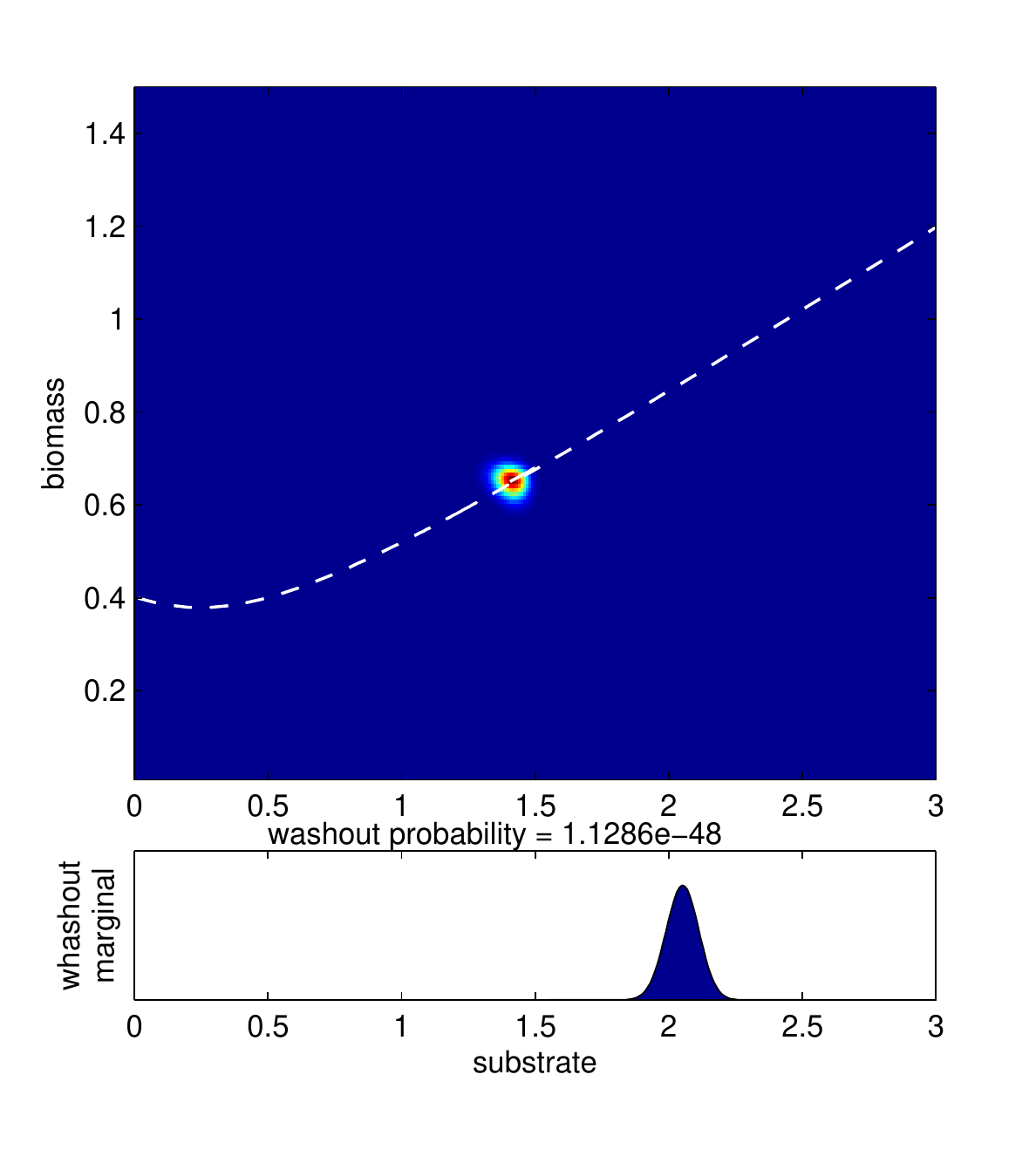}
&
\includegraphics[width=4.7cm,keepaspectratio]{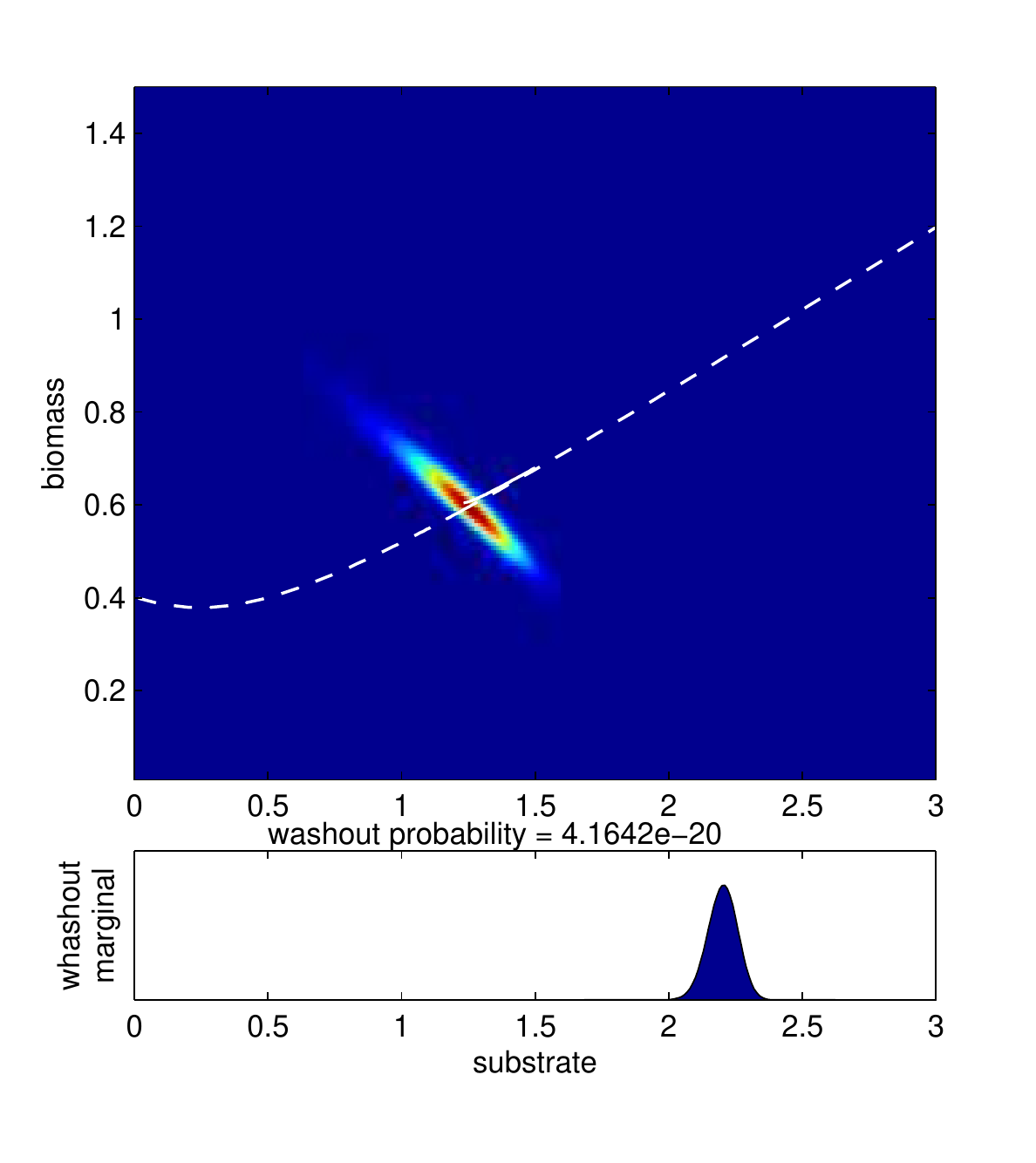}
\\[-1.3em]
$t=0$
&
$t=4$
&
$t=24$
\\ 
\includegraphics[width=4.7cm,keepaspectratio]{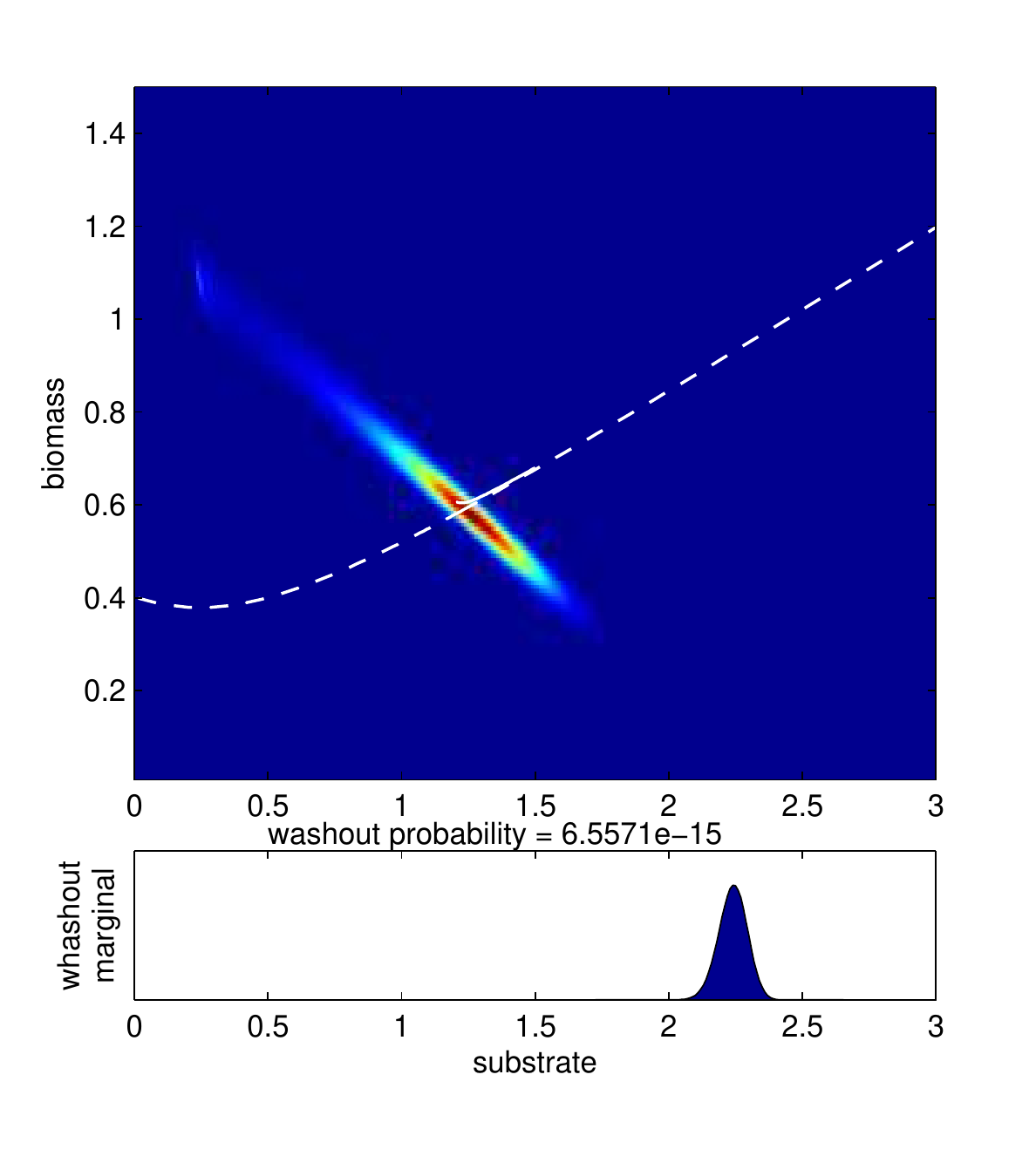}
&
\includegraphics[width=4.7cm,keepaspectratio]{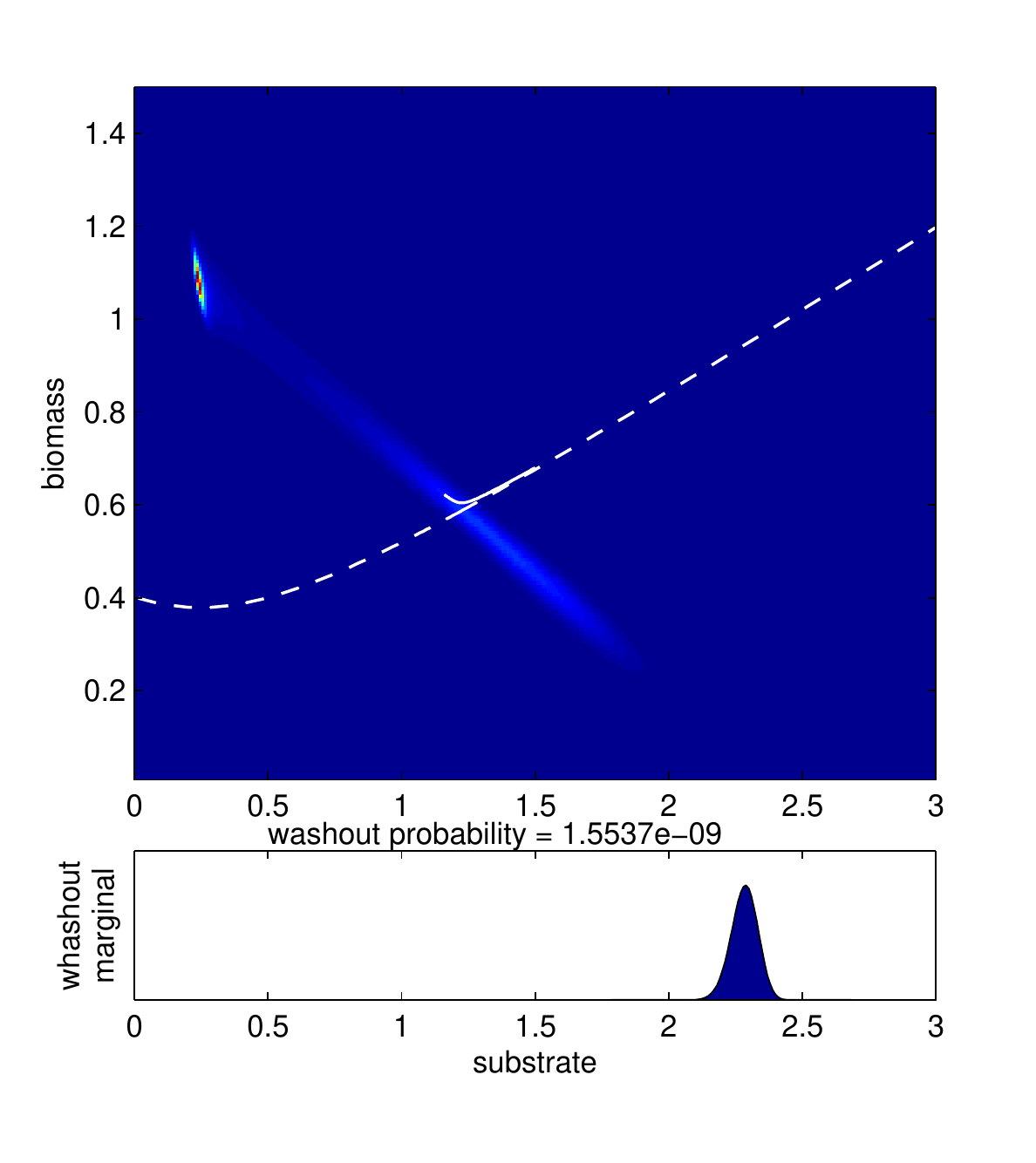}
&
\includegraphics[width=4.7cm,keepaspectratio]{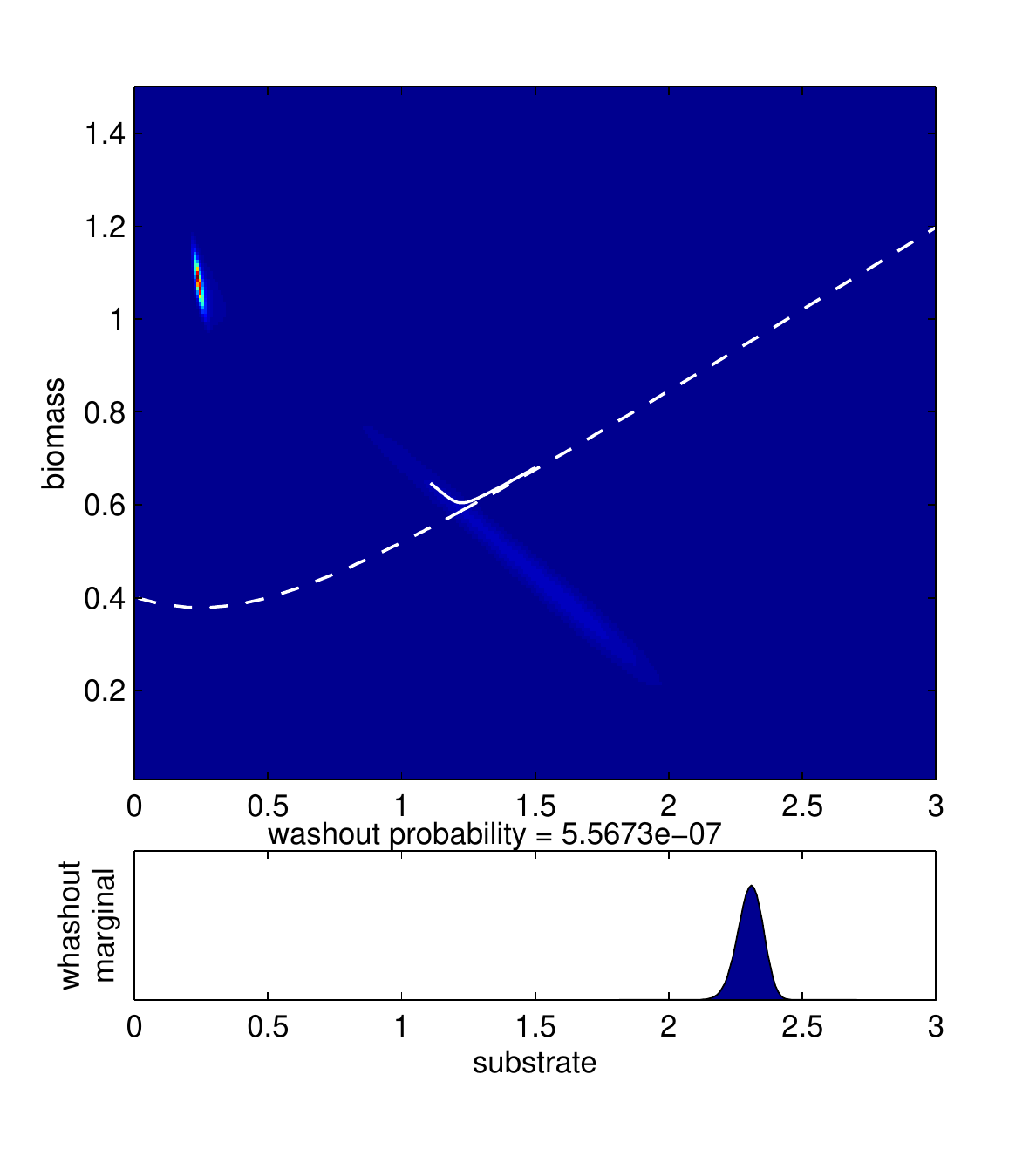}
\\[-1.3em]
$t=32$
&
$t=44$
&
$t=52$
\\ 
\includegraphics[width=4.7cm,keepaspectratio]{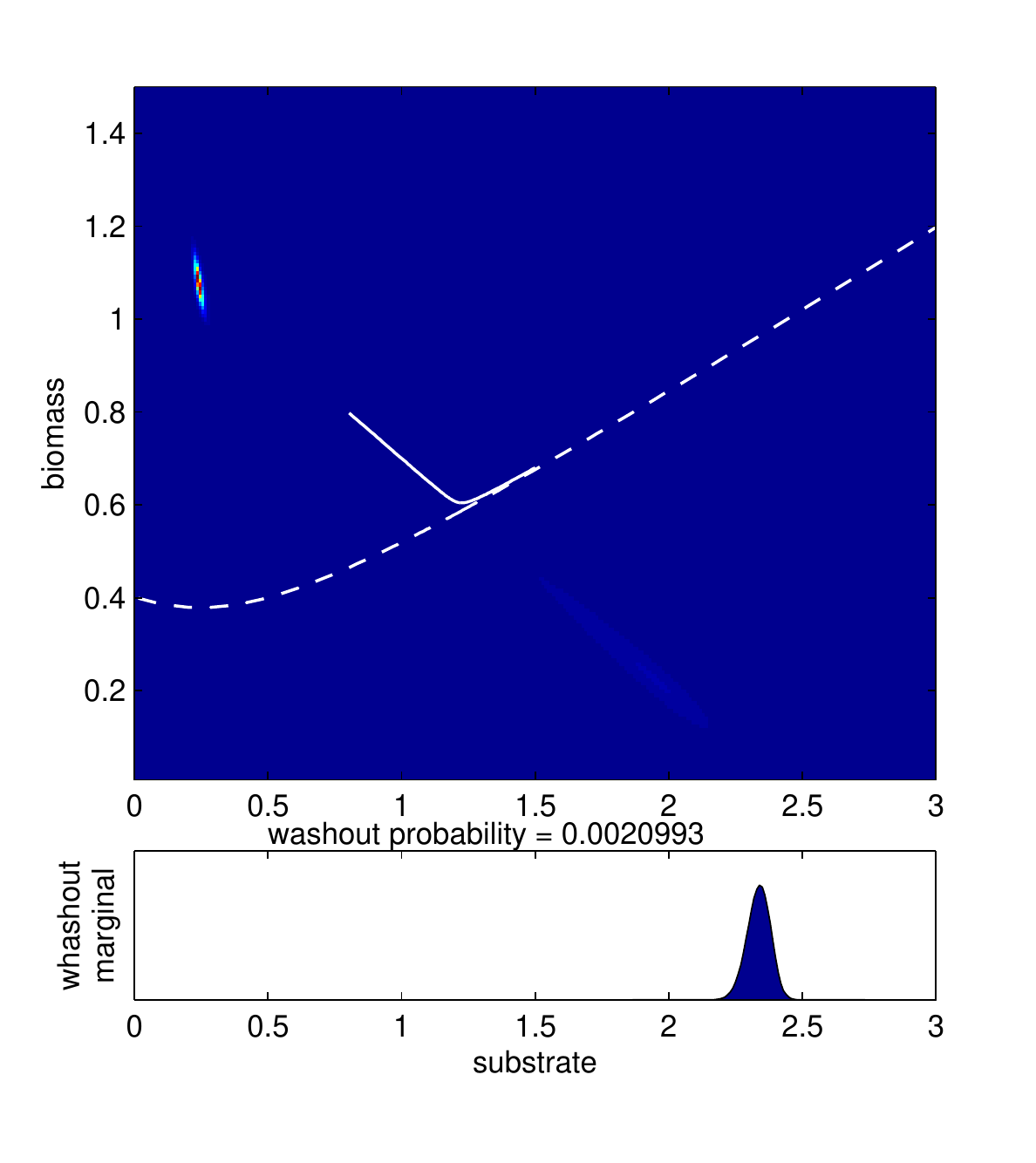}
&
\includegraphics[width=4.7cm,keepaspectratio]{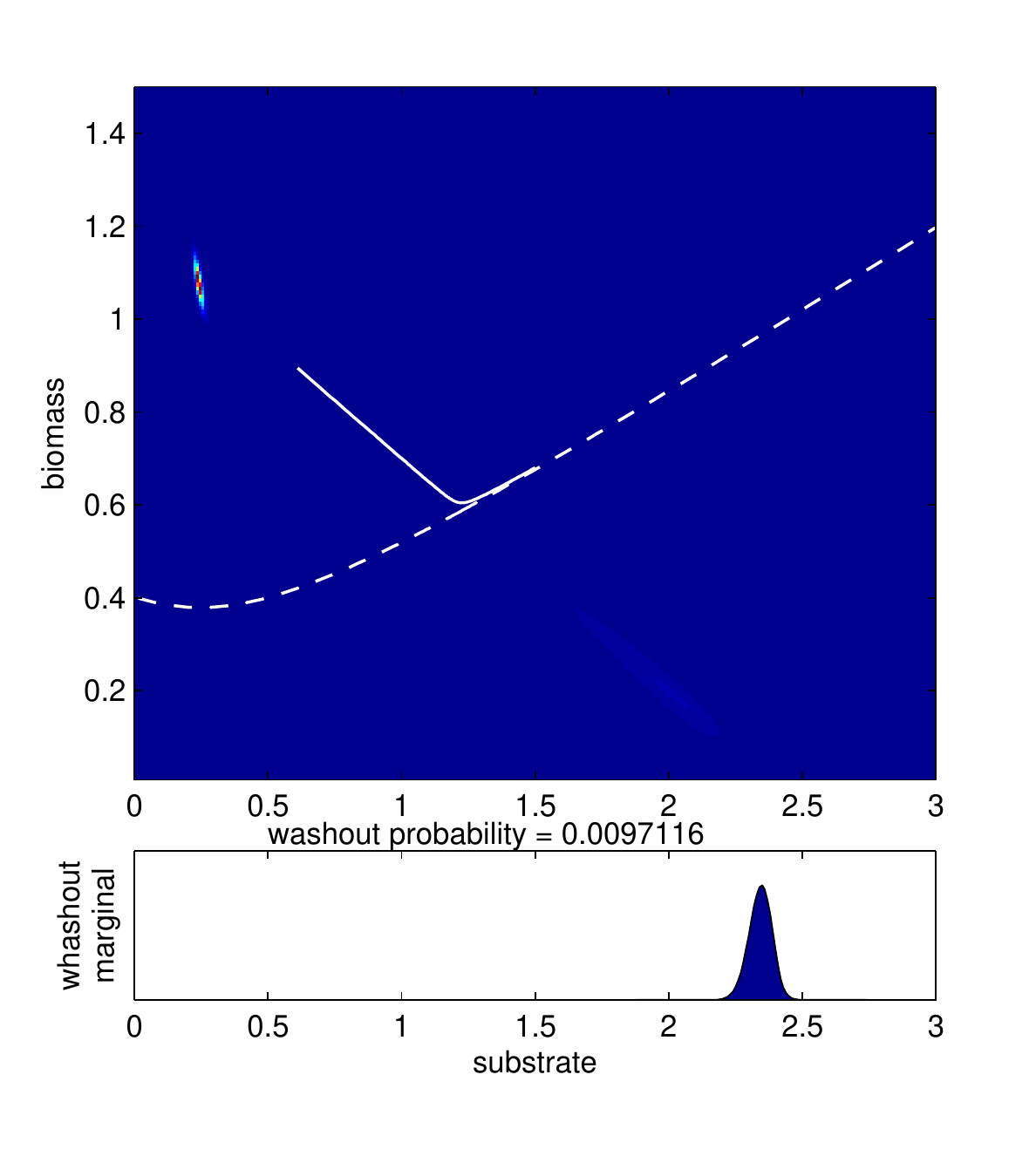}
&
\includegraphics[width=4.7cm,keepaspectratio]{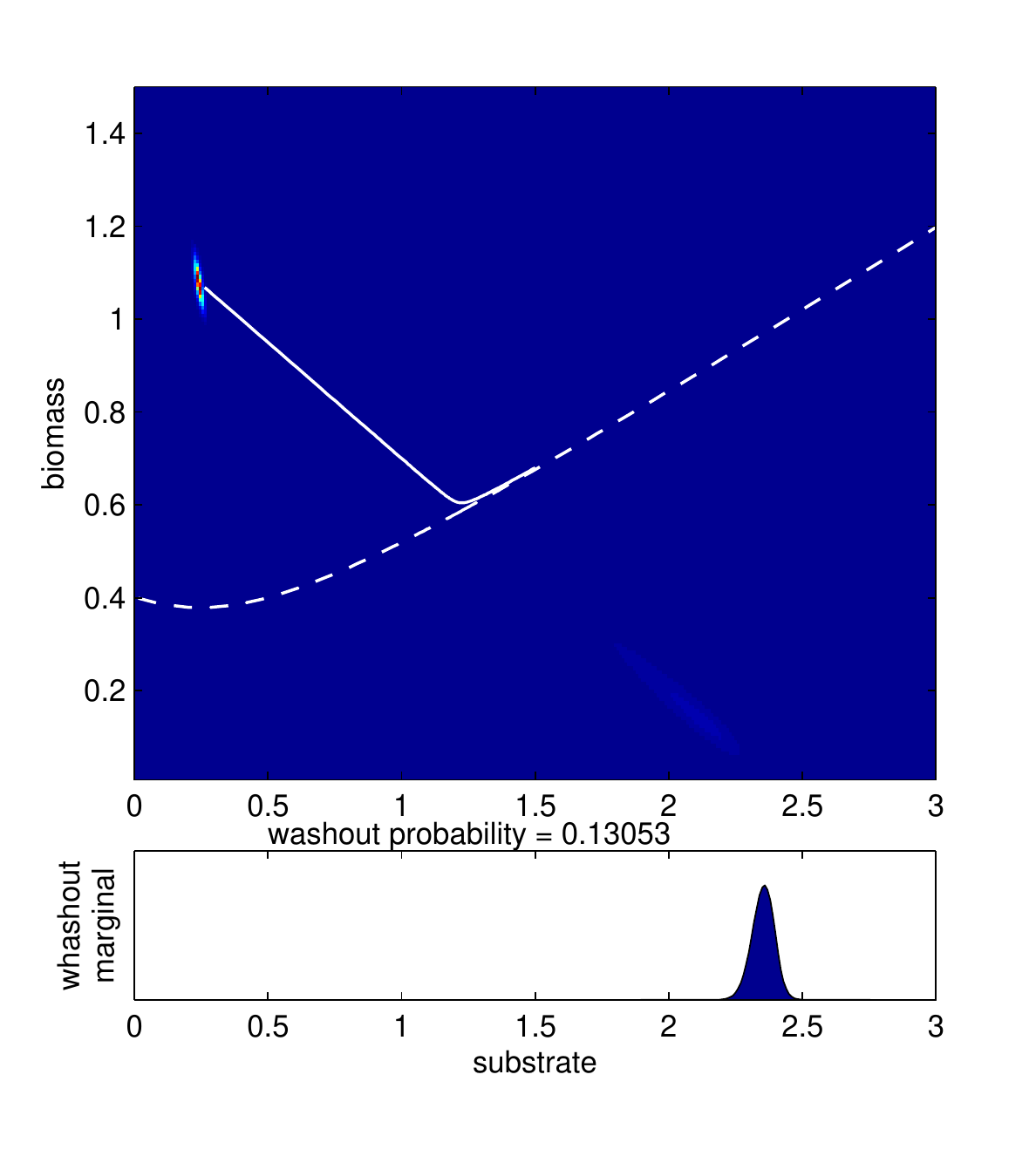}
\\[-1.3em]
$t=68$
&
$t=72$
&
$t=80$
\end{tabular}
\end{center}
}
\caption{Evolution of the distribution law of $X_{t}$: for each time $t$, the density $p_{t}(s,b)$ together with the washout density $q_{t}(s)$; the dashed curve separates the two basins of attraction. The mean of $X_{0}$ is on this curve.
See comments in the text.}
\label{fig.haldane.test}
\end{figure}

In Figure \ref{fig.haldane.test} we plot the time evolution of the distribution law of $X_{t}$: for each time $t$, we represent (the approximation of) $(p_{t}(s,b);(s,b)\in(0,\smax)\times (0,\bmax))$ together with (the approximation of) $(q_{t}(s);s\in(0,\smax))$. In this test the mean of $X_{0}$ is on this curve that separates the two basins of attraction (dashed white line): hence part of the mass will be attracted by $(s^*_{1},b^*_{1})$ and the other part will be attracted by the washout $(\Sin,0)$ (see Figure \ref{fig.monod.haldane}).

For $t=0$ we plot all the trajectory $(x(t))_{t\in[0;80]}$ (white line). At the beginning the distribution law starts to ``stretch'' between the two attractors ($t=24$). At $t=32$, part of the mass is already on the point $(s^*_{1},b^*_{1})$. Note that at this instant $p_{t}(s,b)$ is bimodal and $x(t)$ is a good approximation of $\E(X_{t})$, but it is a poor statistics for $X_{t}$. At the final time $t=80$, the deterministic trajectory $x(t)$ reaches the equilibrium point $(s^*_{1},b^*_{1})$ and $13\%$ of the mass has been trapped by the washout absorbing boundary and some mass is still in the washout basin and will be trapped by the boundary ``$b=0$''.

\appendix

\section{General finite difference scheme for $n$-dimensional diffusion processes}
\label{sec.DF.appendix}

Let $X_{t}$ be the following diffusion process:
\begin{align*}
  \rmd X_{t}
  &=
  b(X_{t})\,\rmd t
  +
  \sigma(X_{t})\,\rmd W_t
\end{align*}
where $X_{t}$ takes values in $\R^n$, $b:\,\R^n\mapsto\R^n$, $\sigma:\,\R^n\mapsto\R^{n\times m}$, and $W_{t}$ is a standard Brownian motion with values in $\R^m$. Let $a=\sigma\,\sigma^*:\,\R^n\mapsto\R^{n\times n}$. The coefficients are supposed 
to be locally Lipschitz and at most of linear growth.

The probability density function $p(t,x)$ of $X_{t}$ is solution of the following Fokker-Planck equation:
\begin{align}
\label{eq.fp}
  \frac{\partial}{\partial t}\,p(t,x)
  &=
  \LL^* p(t,x)
\end{align}
where $\LL$ is the infinitesimal generator defined by:
\begin{align*}
  \LL \phi(x)
  &\eqdef
  \sum_{i=1}^n f_i(x)\,\phi'_{x_{i}}(x)
  +
  \frac12\,\sum_{i,j=1}^n a_{ij}(x)\,\phi''_{x_{i}x_{j}}(x)\,.
\end{align*}
We consider finite difference schemes based on the following stencil (for the components $(x_{i},x_{j})$):
\begin{center}
\includegraphics[width=4.5cm,keepaspectratio]{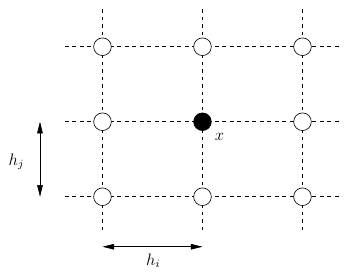}
\end{center}
We use the following up-wind scheme \cite{kushner1977a}:
\begin{align*}
   f_i(x)\,\phi'_{x_{i}}(x)
   &\simeq
   \begin{cases}
      f_i(x)\,\frac{\phi(x+h_{i}\,e_{i})-\phi(x)}{h_{i}}\,,
      &\textrm{ if } f_i(x)\geq 0\,,
   \\
      f_i(x)\,\frac{\phi(x)-\phi(x-h_{i}\,e_{i})}{h_{i}}\,,
      &\textrm{ if } f_i(x)< 0\,,
   \end{cases}
\\[1em]
   a_{ii}(x)\,\phi''_{x_{i}^2}(x)
   &\simeq
   \textstyle
   a_{ii}(x)\,\frac{\phi(x+h_{i}\,e_{i})-2\,\phi(x)+\phi(x-h_{i}\,e_{i})}{h_{i}^2}\,,
\\[1em]
   a_{ij}(x)\,\phi''_{x_{i}x_{j}}(x)
   &\simeq
   \begin{cases}
   \textstyle
      a_{ij}(x)\,\frac{1}{2\,h_{i}}\Big[
         \frac{\phi(x+h_{i}\,e_{i}+h_{j}\,e_{j})-\phi(x+h_{i}\,e_{i})}{h_{j}}
         -
         \frac{\phi(x+h_{j}\,e_{j})-\phi(x)}{h_{j}}
   \\
   \qquad\qquad\qquad
   \textstyle
         +
         \frac{\phi(x)-\phi(x-h_{j}\,e_{j})}{h_{j}}
         -
         \frac{\phi(x-h_{i}\,e_{i})-\phi(x-h_{i}\,e_{i}-h_{j}\,e_{j})}{h_{j}}
      \Big]\,,
   \\
   \qquad\qquad\qquad\qquad\qquad\qquad\qquad\qquad\qquad
   \textrm{ if }a_{ij}(x)\geq 0\,,
\\[1em]
   \textstyle
      a_{ij}(x)\,\frac{1}{2\,h_{i}}\Big[
         \frac{\phi(x+h_i\,e_i)-\phi(x+h_i\,e_i-h_j\,e_j)}{h_{j}}
         -
         \frac{\phi(x)-\phi(x-h_j\,e_j)}{h_{j}}
   \\
   \qquad\qquad\qquad
   \textstyle
         +
         \frac{\phi(x+h_j\,e_j)-\phi(x)}{h_{j}}
         -
         \frac{\phi(x-h_i\,e_i+h_j\,e_j)-\phi(x-h_i\,e_i)}{h_{j}}
      \Big]\,,
   \\
   \qquad\qquad\qquad\qquad\qquad\qquad\qquad\qquad\qquad
   \textrm{ if }a_{ij}(x)< 0\,,
   \end{cases}
\end{align*}
for $i,j=1,\dots,n$, $i\neq j$.
The last non-diagonal second order schemes correspond to the following diagrams:
\begin{center}
\includegraphics[width=8.5cm,keepaspectratio]{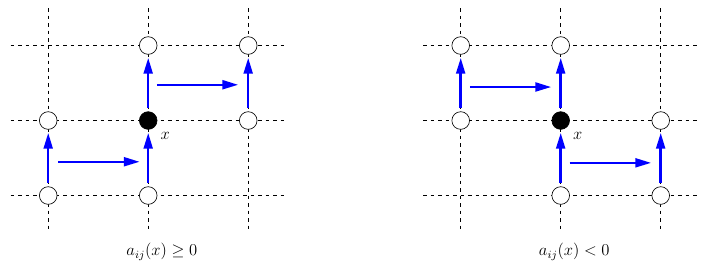}
\end{center}

With notation $f^+(x)=\max(f(x),0)$ and $f^-(x)=\max(-f(x),0)$, we get the following approximation:
{\small
\begin{align*}
   &\LL_h \phi(x)
   =
   \textstyle
   \sum_{i }  f_i(x)\,\phi'_{x_{i}}(x)
   +
   \frac12\,\sum_{i,j }  a_{i,j}(x)\,\phi''_{x_{i}x_{j}}(x)
\\[1em]
   &\quad=
   \textstyle
   \sum_{i }  \{f_i^+(x)-f_i^-(x)\}\,\phi'_{x_{i}}(x)
   +
   \frac12\,\sum_{i }  a_{ii}(x)\,\phi''_{x_{i}^2}(x)
   +
   \frac12\,\sum_{i,j;i\neq j} \{a_{ij}^+(x)-a_{ij}^-(x)\}\,\phi''_{x_{i}x_{j}}(x)
\\[1em] 
   &
   \quad 
   \simeq
   \sum_{i}
   \Big\{
       \textstyle
       \frac{f_i^+(x)}{h_{i}}\,[\phi(x+h_{i}\,e_{i})-\phi(x)]
       -
       \frac{f_i^-(x)}{h_{i}}\,[\phi(x)-\phi(x-h_{i}\,e_{i})]
    \Big\}
\\
   &
   \qquad\qquad 
   +\sum_{i} 
      \textstyle
      \frac{a_{ii}(x)}{2\,h_{i}^2}\,
      [\phi(x+h_{i}\,e_{i})-2\,\phi(x)+\phi(x-h_{i}\,e_{i})]
\\
   &
   \qquad \qquad
   +\frac12\,\sum_{i,j;i\neq j} 
   \textstyle
   \Big\{
   \frac{a_{ij}^+(x)}{2\,h_{i}\,h_{j}}\,
   \Big(
         [\phi(x+h_{i}\,e_{i}+h_{j}\,e_{j})-\phi(x+h_{i}\,e_{i})]
         -
         [\phi(x+h_{j}\,e_{j})-\phi(x)]
   \\[-1em]
   &
   \qquad\qquad\qquad\qquad\qquad\qquad
         \textstyle
         +
         [\phi(x)-\phi(x-h_{j}\,e_{j})]
         -
         [\phi(x-h_{i}\,e_{i})-\phi(x-h_{i}\,e_{i}-h_{j}\,e_{j})]
      \Big)
\\
   &\textstyle
   \qquad\qquad\qquad
      -
      \frac{a_{ij}^-(x)}{2\,h_{i}\,h_{j}}
      \Big(
         [\phi(x+h_i\,e_i)-\phi(x+h_i\,e_i-h_j\,e_j)]
         -
         [\phi(x)-\phi(x-h_j\,e_j)]
   \\
   &
   \qquad\qquad\qquad\qquad\qquad\qquad
         \textstyle
         +
         [\phi(x+h_j\,e_j)-\phi(x)]
         -
         [\phi(x-h_i\,e_i+h_j\,e_j)-\phi(x-h_i\,e_i)]
      \Big)
      \Big\}
\\[1em] 
&=  
   \phi(x)\;
   \textstyle
   \Big\{
     -\sum_{i}\frac{|f_{i}(x)|}{h_{i}}
     -\sum_{i}\frac{a_{ii}(x)}{h_{i}^2}
     +\sum_{i,j;i\neq j}\frac{|a_{ij}(x)|}{2\,h_{i}\,h_{j}}
   \Big\}
\\
&\qquad
   +
   \textstyle\sum_{i}
   \phi(x+h_{i}\,e_{i})\;
   \textstyle
   \Big\{
       \frac{f_{i}^+(x)}{h_{i}}
     + \frac{a_{ii}(x)}{2\,h_{i}^2}
     -\sum_{j;j\neq i}\frac{|a_{ij}(x)|}{4\,h_{i}\,h_{j}}
   \Big\}
   +
   \sum_{j}
   \textstyle
    \sum_{i;i\neq j}
       \phi(x+h_{j}\,e_{j})\;\frac{|a_{ij}(x)|}{4\,h_{i}\,h_{j}}
\\
&\qquad
   +
   \textstyle\sum_{i}
   \phi(x-h_{i}\,e_{i})\;
   \textstyle
   \Big\{
       \frac{f_{i}^-(x)}{h_{i}}
     + \frac{a_{ii}(x)}{2\,h_{i}^2}
     -\sum_{j;j\neq i}\frac{|a_{ij}(x)|}{4\,h_{i}\,h_{j}}
   \Big\}
   +
   \sum_{j}
   \textstyle
    \sum_{i;i\neq j}
       \phi(x-h_{j}\,e_{j})\;\frac{|a_{ij}(x)|}{4\,h_{i}\,h_{j}}
\\
   &
   \qquad
   +\textstyle\sum_{i,j;i\neq j} 
   \textstyle
   \Big\{
   \frac{a_{ij}^+(x)}{2\,h_{i}\,h_{j}}\,
   [\phi(x+h_{i}\,e_{i}+h_{j}\,e_{j})+\phi(x-h_{i}\,e_{i}-h_{j}\,e_{j})]
\\
   &\textstyle
   \qquad\qquad\qquad\qquad\qquad
      +
      \frac{a_{ij}^-(x)}{2\,h_{i}\,h_{j}}
      [\phi(x+h_i\,e_i-h_j\,e_j)+\phi(x-h_i\,e_i+h_j\,e_j)]
      \Big\}
\end{align*}
}
\noindent the symmetry $a_{ij}=a_{ji}$ leads to
{\small
\begin{align*}
   \LL_h \phi(x)
   &=
   \textstyle
   \phi(x)\;
   \textstyle
   \Big\{
     -\sum_{i}\frac{|f_{i}(x)|}{h_{i}}
     -\sum_{i}\frac{a_{ii}(x)}{h_{i}^2}
     +\sum_{i,j;i\neq j}\frac{|a_{ij}(x)|}{2\,h_{i}\,h_{j}}
   \Big\}
\\
&\qquad
   +
   \textstyle\sum_{i}
   \phi(x+h_{i}\,e_{i})\;
   \textstyle
   \Big\{
       \frac{f_{i}^+(x)}{h_{i}}
     + \frac{a_{ii}(x)}{2\,h_{i}^2}
     -\sum_{j;j\neq i}\frac{|a_{ij}(x)|}{2\,h_{i}\,h_{j}}
   \Big\}
\\
&\qquad
   +
   \textstyle\sum_{i}
   \phi(x-h_{i}\,e_{i})\;
   \textstyle
   \Big\{
       \frac{f_{i}^-(x)}{h_{i}}
     + \frac{a_{ii}(x)}{2\,h_{i}^2}
     -\sum_{j;j\neq i}\frac{|a_{ij}(x)|}{2\,h_{i}\,h_{j}}
   \Big\}
\\
   &
   \qquad
   +\textstyle\sum_{i,j;i\neq j} 
   \textstyle
   \Big\{
   \frac{a_{ij}^+(x)}{2\,h_{i}\,h_{j}}\,
   [\phi(x+h_{i}\,e_{i}+h_{j}\,e_{j})+\phi(x-h_{i}\,e_{i}-h_{j}\,e_{j})]
\\
   &\textstyle
   \qquad\qquad\qquad\qquad\qquad
      +
      \frac{a_{ij}^-(x)}{2\,h_{i}\,h_{j}}
      [\phi(x+h_i\,e_i-h_j\,e_j)+\phi(x-h_i\,e_i+h_j\,e_j)]
      \Big\}
\end{align*}
}
We get the following approximation of the infinitesimal generator:
\[
  \LL \phi(x)
  \simeq
  \LL_{h} \phi(x)
  =
  \sum_{y\in G_{h}}\LL_h(x,y)\,\phi(y)
\]
for $x\in G_{h}$ where $G_{h}=\{x=(k_{1}\,h_{1},\dots,k_{n}\,h_{n})\,;\,k_{i}=0,\dots,N_{i},\, i=1,\dots,n\}$ and
\begin{align*}
\left\{
\begin{array}{rll}
  \LL_{h} (x,x)
  &=
  -\sum_{i=1}^n  \frac{|f_{i}(x)|}{h_{i}}
  -\sum_{i=1}^n \Big\{
    \frac{a_{ii}(x)}{h_{i}^2} - \sum_{j\neq i}\frac{|a_{ij}(x)|}{2\,h_{i}\,h_{j}}
  \Big\}
  \,,
\\
  \LL_{h} (x,x\pm h_i\,e_i)
  &=
  \frac{f_{i}^\pm(x)}{h_{i}}
  +
  \frac{a_{ii}(x)}{2\,h_{i}^2} - \sum_{j;j\neq i}\frac{|a_{ij}(x)|}{2\,h_{i}\,h_{j}} 
     \,,
\\
  \LL_{h} (x,x+ h_i\,e_i+ h_j\,e_j)
  &=
  \LL_{h} (x,x- h_i\,e_i- h_j\,e_j)
  =
  \frac{a_{ij}^+(x)}{2\,h_{i}\,h_{j}} 
  \qquad\textrm{ for }i\neq j\,,
\\
  \LL_{h} (x,x+ h_i\,e_i- h_j\,e_j)
  &=
  \LL_{h} (x,x- h_i\,e_i+ h_j\,e_j)
  =
  \frac{a_{ij}^-(x)}{2\,h_{i}\,h_{j}} 
  \qquad\textrm{ for }i\neq j\,,
\\[0.7em]
  \LL_{h} (x,y)
  &=
  0 \qquad\textrm{otherwise.}
\end{array}
\right.
\end{align*}

\section{Boundary conditions for the finite difference approximation}
\label{sec.bc}

For the boundary points $G_{h}\setminus \mathring G_{h}$ of the grid, we use the following schemes:
\begin{itemize}
\item For $x\in\{(s,b)\in G_{h}\,;\,s=0,b\in(0,\bmax)\}$
\begin{align*}
\left\{
\begin{array}{rll}
  \LL_{h} (x,x)
  &=
  - \frac{|f_1(x)|}{h_1} - \frac{|f_2(x)|}{h_2} - \frac{\sigma^2_1(x)}{h_1^2}
  - \frac{\sigma^2_2(x)}{h_2^2}
  \,,
\\
  \LL_{h} (x,x+ h_1\,e_1)
  &=
  \frac{f_{1}^+(x)}{h_{1}}
  +
  \frac{\sigma^2_{1}(x)}{2\,h_{1}^2} 
     \,,
\\
  \LL_{h} (x,x- h_1\,e_1)
  &=
  \frac{f_{1}^-(x)}{h_{1}}
  +
  \frac{\sigma^2_{1}(x)}{2\,h_{1}^2} 
  =
  0
  \textrm{ because }f_{1}(0,b)=D\,\Sin\,,\ \sigma_{1}(0,b)=0\,,
\\
  \LL_{h} (x,x\pm h_2\,e_2)
  &=
  \frac{f_{2}^\pm(x)}{h_{2}}
  +
  \frac{\sigma^2_{2}(x)}{2\,h_{2}^2} 
     \,,
     \textrm{ note that }f_{2}(0,b)=-D\,b<0\,,
\\[0.7em]
  \LL_{h} (x,y)
  &=
  0 \qquad\textrm{otherwise.}
\end{array}
\right.
\end{align*}
\item For $x\in\{(s,b)\in G_{h}\,;\,s=\smax,b\in(0,\bmax)\}$
\begin{align*}
\left\{
\begin{array}{rll}
  \LL_{h} (x,x)
  &=
  - \frac{|f_1(x)|}{h_1} - \frac{|f_2(x)|}{h_2} - \frac{\sigma^2_1(x)}{h_1^2}
  - \frac{\sigma^2_2(x)}{h_2^2}
  \,,
\\
  \LL_{h} (x,x+ h_1\,e_1)
  &=
  0\textrm{ (set artificially to $0$)} 
     \,,
\\
  \LL_{h} (x,x- h_1\,e_1)
  &=
  \frac{|f_{1}(x)|}{h_{1}}
  +
  \frac{\sigma^2_{1}(x)}{h_{1}^2} 
     \,,
\\
  \LL_{h} (x,x\pm h_2\,e_2)
  &=
  \frac{f_{2}^\pm(x)}{h_{2}}
  +
  \frac{\sigma^2_{2}(x)}{2\,h_{2}^2} 
     \,,
\\[0.7em]
  \LL_{h} (x,y)
  &=
  0 \qquad\textrm{otherwise.}
\end{array}
\right.
\end{align*}
\item For $x\in\{(s,b)\in G_{h}\,;\,s\in(0,\smax),b=0\}$
\begin{align*}
\left\{
\begin{array}{rll}
  \LL_{h} (x,x)
  &=
  - \frac{|f_1(x)|}{h_1}  - \frac{\sigma^2_1(x)}{h_1^2}
  \,,
\\
  \LL_{h} (x,x\pm h_1\,e_1)
  &=
  \frac{f_{1}^\pm(x)}{h_{1}}
  +
  \frac{\sigma^2_{1}(x)}{2\,h_{1}^2} 
     \,,
\\
  \LL_{h} (x,x\pm h_2\,e_2)
  &=
  \frac{f_{2}^\pm(x)}{h_{2}}
  +
  \frac{\sigma^2_{2}(x)}{2\,h_{2}^2} 
  =0
  \textrm{ because $f_{2}(s,0)=\sigma_{2}(s,0)=0$}
     \,,
\\[0.7em]
  \LL_{h} (x,y)
  &=
  0 \qquad\textrm{otherwise.}
\end{array}
\right.
\end{align*}
\item For $x\in\{(s,b)\in G_{h}\,;\,s\in(0,\smax),b=\bmax\}$
\begin{align*}
\left\{
\begin{array}{rll}
  \LL_{h} (x,x)
  &=
  - \frac{|f_1(x)|}{h_1} - \frac{|f_2(x)|}{h_2} - \frac{\sigma^2_1(x)}{h_1^2}
  - \frac{\sigma^2_2(x)}{h_2^2}
  \,,
\\
  \LL_{h} (x,x\pm h_1\,e_1)
  &=
  \frac{f_{1}^\pm(x)}{h_{1}}
  +
  \frac{\sigma^2_{1}(x)}{2\,h_{1}^2} 
     \,,
\\
  \LL_{h} (x,x+ h_2\,e_2)
  &=
  0 \textrm{ (set artificially to $0$)}
     \,,
\\
  \LL_{h} (x,x- h_2\,e_2)
  &=
  \frac{|f_{2}(x)|}{h_{2}}
  +
  \frac{\sigma^2_{2}(x)}{h_{2}^2} 
     \,,
\\[0.7em]
  \LL_{h} (x,y)
  &=
  0 \qquad\textrm{otherwise.}
\end{array}
\right.
\end{align*}
\item For $x=(0,0)$
\begin{align*}
\left\{
\begin{array}{rll}
  \LL_{h} (x,x)
  &=
  - \frac{|f_1(x)|}{h_1}  
  \,,
\\
  \LL_{h} (x,x+ h_1\,e_1)
  &=
  \frac{f_{1}^+(x)}{h_{1}}
     \,,
\\
  \LL_{h} (x,x- h_1\,e_1)
  &=
  0
  \textrm{ because }f_{1}(0,0)=D\,\Sin\,,\ \sigma_{1}(0,0)=0\,,
     \,,
\\
  \LL_{h} (x,x\pm h_2\,e_2)
  &=
  0 \textrm{ because }f_{2}(0,0)=\sigma_{2}(0,0)=0 
     \,,
\\[0.7em]
  \LL_{h} (x,y)
  &=
  0 \qquad\textrm{otherwise.}
\end{array}
\right.
\end{align*}
\item For $x=(\smax,0)$
\begin{align*}
\left\{
\begin{array}{rll}
  \LL_{h} (x,x)
  &=
  - \frac{|f_1(x)|}{h_1} - \frac{\sigma^2_1(x)}{h_1^2}
  \,,
\\
  \LL_{h} (x,x+ h_1\,e_1)
  &=
  0\textrm{ (set artificially to $0$)}
     \,,
\\
  \LL_{h} (x,x- h_1\,e_1)
  &=
  \frac{|f_{1}(x)|}{h_{1}}
  +
  \frac{\sigma^2_{1}(x)}{h_{1}^2} 
     \,,
\\
  \LL_{h} (x,x\pm h_2\,e_2)
  &=
  0 \textrm{ because }f_{2}(\smax,0)=\sigma_{2}(\smax,0)=0 
     \,,
\\[0.7em]
  \LL_{h} (x,y)
  &=
  0 \qquad\textrm{otherwise.}
\end{array}
\right.
\end{align*}
\item For $x=(0,\bmax)$
\begin{align*}
\left\{
\begin{array}{rll}
  \LL_{h} (x,x)
  &=
  - \frac{|f_1(x)|}{h_1} - \frac{|f_2(x)|}{h_2} 
  - \frac{\sigma^2_2(x)}{h_2^2}
  \,,
\\ 
  \LL_{h} (x,x+ h_1\,e_1)
  &=
  \frac{f_{1}^+(x)}{h_{1}}
     \,,
\\ 
  \LL_{h} (x,x- h_1\,e_1)
  &=
  0\textrm{ because }f_{1}(0,\bmax)=D\,\Sin\,,\ \sigma_{1}(0,\bmax)=0
     \,,
\\
  \LL_{h} (x,x+ h_2\,e_2)
  &=
  0\textrm{ (set artificially to $0$)}
     \,,
\\
  \LL_{h} (x,x- h_2\,e_2)
  &=
  \frac{|f_{2}(x)|}{h_{2}}
  +
  \frac{\sigma^2_{2}(x)}{h_{2}^2} 
     \,,
\\[0.7em]
  \LL_{h} (x,y)
  &=
  0 \qquad\textrm{otherwise.}
\end{array}
\right.
\end{align*}
\item For $x=(\smax,\bmax)$
\begin{align*}
\left\{
\begin{array}{rll}
  \LL_{h} (x,x)
  &=
  - \frac{|f_1(x)|}{h_1} - \frac{|f_2(x)|}{h_2} - \frac{\sigma^2_1(x)}{h_1^2}
  - \frac{\sigma^2_2(x)}{h_2^2}
  \,,
\\
  \LL_{h} (x,x+ h_1\,e_1)
  &=
  0\textrm{ (set artificially to $0$)}
     \,,
\\
  \LL_{h} (x,x- h_1\,e_1)
  &=
  \frac{|f_{1}(x)|}{h_{1}}
  +
  \frac{\sigma^2_{1}(x)}{h_{1}^2} 
     \,,
\\
  \LL_{h} (x,x+ h_2\,e_2)
  &=
  0\textrm{ (set artificially to $0$)}
     \,,
\\
  \LL_{h} (x,x- h_2\,e_2)
  &=
  \frac{|f_{2}(x)|}{h_{2}}
  +
  \frac{\sigma^2_{2}(x)}{h_{2}^2} 
     \,,
\\[0.7em]
  \LL_{h} (x,y)
  &=
  0 \qquad\textrm{otherwise.}
\end{array}
\right.
\end{align*}
\end{itemize}



\section*{Acknowledgements}

The work was partially supported by the French National Research Agency (ANR) within the SYSCOMM project ANR-09- SYSC-003.



\addcontentsline{toc}{section}{Reference}
\bibliographystyle{plain}

\begin{thebibliography}{}

\end{thebibliography}


\begin{thebibliography}{10}

\bibitem{brezis2010a}
H.~Br{\'e}zis.
\newblock {\em Functional Analysis, {Sobolev} Spaces and Partial Differential
  Equations}.
\newblock Springer, 2010.

\bibitem{campillo2011chemostat}
Fabien Campillo, Marc Joannides, and Ir{\`e}ne Larramendy-Valverde.
\newblock Stochastic modeling of the chemostat.
\newblock {\em Ecological Modelling}, 222(15):2676--2689, 2011.

\bibitem{campillo2012analysischemostat}
Fabien Campillo, Marc Joannides, and Ir{\`e}ne Larramendy-Valverde.
\newblock Analysis of the stochastic chemostat.
\newblock In preparation, 2012.

\bibitem{grasman1999asymptotic}
J.~Grasman and O.A. Herwaarden.
\newblock {\em Asymptotic methods for the Fokker-Planck equation and the exit
  problem in applications}.
\newblock Springer, 1999.

\bibitem{grasman2005a}
Johan Grasman and Maarten {De Gee}.
\newblock Breakdown of a chemostat exposed to stochastic noise volume.
\newblock {\em Journal of Engineering Mathematics}, 53(3):291--300, 2005.

\bibitem{ikeda1981a}
Nobuyuki Ikeda and Shinzo Watanabe.
\newblock {\em Stochastic Differential Equations and Diffusion Processes}.
\newblock North--Holland/Kodansha, Amsterdam, 1981.

\bibitem{imhof2005a}
Lorens Imhof and Sebastian Walcher.
\newblock Exclusion and persistence in deterministic and stochastic chemostat
  models.
\newblock {\em Journal of Differential Equations}, 217(1):26--53, 2005.

\bibitem{kushner1977a}
Harold~J. Kushner.
\newblock {\em Probability Methods for Approximations in Stochastic Control and
  for Elliptic Equations}, volume 129 of {\em Mathematics in Science and
  Engineering}.
\newblock Academic Press, New York, 1977.

\bibitem{kushner1990a}
Harold~J. Kushner.
\newblock Numerical methods for stochastic control problems in continuous time.
\newblock {\em SIAM J. Control Optim.}, 28(5):999--1048, 1990.

\bibitem{lamberton1996a}
Damien Lamberton and Bernard Lapeyre.
\newblock {\em Introduction to Stochastic Calculus Applied to Finance}.
\newblock {Chapman \&\ Hall/CRC}, 1996.

\bibitem{schuss2010a}
Zeev Schuss.
\newblock {\em Theory and Applications of Stochastic Processes, An Analytical
  Approach}.
\newblock Springer, 2010.

\bibitem{smith1995a}
Hal~L. Smith and Paul~E. Waltman.
\newblock {\em The Theory of the Chemostat: Dynamics of Microbial Competition}.
\newblock {Cambridge University Press}, 1995.

\end{thebibliography}


\end{document}